\documentclass[3p,12pt]{elsarticle}
\usepackage{mathrsfs}
\usepackage{amsmath}
\usepackage{comment}
\usepackage{wrapfig}
\usepackage{graphicx}
\usepackage{amsfonts}
\usepackage{lineno}
\usepackage{epsfig}
\usepackage{amssymb}
\usepackage{amsthm}
\usepackage[colorlinks]{hyperref}

\newtheorem{thm}{Theorem}[section]
\newtheorem{lemma}[thm]{Lemma}
\newtheorem{cor}[thm]{Corollary}
\newtheorem{prop}[thm]{Proposition}
\newtheorem{Def}[thm]{Definition}
\newtheorem{remark}[thm]{Remark}

\def\H{\mathcal{H}}

\def\G{\mathcal{G}}
\def\L{\mathcal{L}}
\def\K{\mathcal{K}}
\def\R{\mathbb{R}}
\def\N{\mathbb{N}}
\def\Z{\mathbb{Z}}
\def\pt{\partial_t}

\def\ptt{\partial_{t}^2}

\def\l({\left(}
\def\r){\right)}
\def\eps{\varepsilon}

\def\al{\alpha}

\def\ga{\gamma}

\def\lm{\lambda}

\def\la{\langle}
\def\ra{\rangle}
\def\td{{\rm d}}
\def\tdx{{\rm d}x}
\def\tdt{{\rm d}t}
\def\sop{\L}
\def\sqsop{\sqrt{\L}}
\def\ssop0{\sqrt{-\Delta}}
\numberwithin{equation}{section}

\journal{}
\begin{document}

\begin{frontmatter}
  \title{Well-posedness and attractors for the sub-quintic weakly damped wave equation with low-regularity forcing on $\R^3$}

  \author{Cuncai Liu$^1$}
  \author{Fengjuan Meng$^1$\corref{cor1}}
  \author{Chunyou Sun$^{2,3}$}

  \address{1. School of Mathematics and Physics, Jiangsu University of Technology, Changzhou, 213001, China\\
    2. School of Mathematics and Statistics, Lanzhou University, Lanzhou, 730000, China\\
    3. School of Mathematics and Statistics, Donghua University, Shanghai, 201620, China}

  \cortext[cor1]{Corresponding author.\\
    E-mails: liucc@jsut.edu.cn (C. Liu), fjmeng@jsut.edu.cn (F. Meng), sunchy@lzu.edu.cn (C. Sun)}
  \begin{abstract}
    We study the weakly damped semilinear wave equation on $\R^3$ with a sub-quintic nonlinearity
    $f(u)$ and a forcing term of low regularity, $g\in H^{-1}(\R^3)$. By absorbing the forcing term
    into a stationary solution $h$ and then expanding the solution around $h$, we reduce the problem
    to a wave equation with potentials. Relying on Strichartz estimates for this equation, obtained
    through the Besov spaces associated with the Schr\"odinger operator, we establish the
    global well-posedness of weak solutions by an iteration argument on the Besov scale. In particular,
    uniqueness holds unconditionally in the energy class for $3<p<5$. Combining a
    refined decomposition of the solution with uniform tail estimates for trajectories starting from
    bounded sets, we then prove that the associated solution semigroup is asymptotically compact
    and, consequently, possesses a compact global attractor in the energy space $H^1(\R^3)\times L^2(\R^3)$.
  \end{abstract}


\end{frontmatter}


\section{Introduction}

The long-time behavior of solutions to the weakly damped semilinear wave equation has been a subject of intense study in the field of infinite-dimensional dynamical systems. Consider the following model defined on the whole space $\R^3$
\begin{eqnarray}\label{eq:problem}
  \left\{
  \begin{array}{lc}
    u_{tt}+\ga u_{t}-\Delta u+m u+f(u)=g,~(t,x)\in \R_{+}\times\R^3, \\
    u(0)=u_0,~~  u_{t}(0)=u_1,   x\in \R^3.
  \end{array}
  \right.
\end{eqnarray}
Here $\ga>0$, $m>0$, and $g\in H^{-1}(\R^3)$ is a time-independent external force. The nonlinearity $f\in C^3(\R)$ satisfies $f(0)=f^{\prime}(0)=0$ together with the following conditions:
\begin{equation*}\label{F1}
  |f^{\prime\prime\prime}(s)|\leq C|s|^{p-3} \text{~as } |s| \ge 1,\tag{F1}
\end{equation*}
where $3<p<5$; for every $s\in \R$,
\begin{equation*}\label{F2}
  f(s)s\ge0;\tag{F2}
\end{equation*}
and, for some constant $\beta\in (0, 1]$,
\begin{equation*}\label{F3}
  f(s)s-\beta F(s)\ge0, \tag{F3}
\end{equation*}
where $F(s)=\int_0^s f(\tau)\td \tau$ is the primitive function of $f$.

We decompose $f$ by a high--low cut-off depending on a large parameter $\Lambda$. Choose $\chi\in C^\infty(\R)$ such that
$$ 0\le\chi\le1,\qquad \chi(s)=0\quad\text{for }|s|\le1/2,\qquad  \chi(s)=1\quad\text{for }|s|\ge1,$$
and set
\begin{equation}\label{eq:fdecomp_def}
  f_1(s):=\chi(s/\Lambda)f(s),\qquad  f_2(s):=(1-\chi(s/\Lambda))f(s).
\end{equation}
Then $f=f_1+f_2$ and
\begin{equation}\label{eq:fdecomp}
  |f_1^{\prime\prime\prime}(s)|\le C_1|s|^{p-3},
  \qquad
  |f_2^\prime(s)|\le C_{\Lambda}:=C_2\Lambda^{p-1}, \qquad s\in\R,
\end{equation}
with $C_1$, $C_2$ independent of $\Lambda$. The value of $\Lambda$ will be
fixed after the stationary solution $h$ has been introduced; see
Lemma~\ref{le_sta_sol} and the paragraph following it.

Throughout the paper, $\la \cdot, \cdot \ra$ denotes the duality pairing between $H^1(\R^3)$ and its dual, and $\|\cdot\|$ the norm of $L^2(\R^3)$.

When the forcing term is sufficiently regular, typically $g\in L^2(\R^3)$, global well-posedness and the existence of a global attractor in the natural energy space
$\H=H^1(\R^3)\times L^2(\R^3)$ are well established. In three-dimensional bounded domains, the cubic growth rate was traditionally regarded as critical: Sobolev embedding
 arguments yielded global well-posedness only for sub-cubic or cubic nonlinearities, and the corresponding attractor theory was likewise restricted to $p\le 3$; see, for instance,
\cite{avb,avb1,Ball04,rt} and the references therein. Advances in harmonic analysis and dispersive equations, particularly Strichartz estimates, have since enabled the treatment of super-cubic and even
quintic nonlinearities. For damped wave equations with supercritical growth, attractor results were obtained in \cite{ef}. For the defocusing quintic wave equation on smooth three-dimensional domains,
Burq, Lebeau, and Planchon established global well-posedness in the energy space using Strichartz estimates \cite{burq08}. Attractor results for damped quintic wave equations in bounded domains were subsequently
established by Kalantarov, Savostianov, and Zelik \cite{zelik16}.

In applications, the force may have low regularity. We consider in particular distributional forces $g\in H^{-1}(\R^3)$, which include some spatially concentrated loads. This low regularity gives rise to
two difficulties.
\begin{enumerate}
  \item \textbf{Energy identity.} Pairing the equation with $u_t$, the standard route to the energy identity, is no longer legitimate: the duality pairing $\la g, u_t \ra$ is undefined for $g \in H^{-1}(\R^3)$ and $u_t \in L^2(\R^3)$.
  \item \textbf{Lack of regularity.} The energy bound gives only $H^1(\R^3)$ regularity for the solution $u(t)$. Although this yields $L^6(\R^3)$ integrability, it does not, in general, provide higher spatial integrability. Consequently, for a super-cubic nonlinearity, $f(u)$ need not belong to $L^2(\R^3)$, and the nonlinear term cannot generally be treated as an $L^2$ source.
\end{enumerate}

This paper treats low-regularity forcing, $g\in H^{-1}(\R^3)$, together with a super-cubic nonlinearity. Wave equations do not provide parabolic smoothing; the energy estimate alone yields $H^1$-regularity, without any general spatial gain. In the super-cubic regime, $f(u)$ then need not lie in $L^2(\R^3)$, while $-\Delta u$ lies only in $H^{-1}(\R^3)$; the space $H^{-1}(\R^3)$ is thus the natural space for $g$. By contrast, in strongly damped wave equations, the extra term $-\Delta u_t$ compensates for the weak regularity of $g$, and both the existence and the asymptotic regularity of global attractors are known in that setting; see \cite{pata,sun10,yang}. However, these methods do not apply to the weakly damped equation \eqref{eq:problem}.

For the weakly damped equation \eqref{eq:problem} with $g\in H^{-1}(\R^3)$, previous results include the following. For cubic nonlinearities with low-regularity forcing, a global attractor
was established on bounded domains in \cite{lcc}. On $\R^3$, \cite{fjm} established a global attractor under its stated assumptions and obtained translational regularity under an additional condition.
For $3<p<4$ on $\R^3$, global well-posedness and attractor results were studied using the notion of translationally regular solutions in \cite{lcc1}. The Dirichlet problem on bounded domains was
 studied in \cite{lcc-jde}; the range considered there is limited by the available Strichartz estimates on bounded domains \cite{sogge09} to $p\le 25/7$.

In the present paper, we establish an unconditional well-posedness theory for weak solutions in the range $3<p<5$. In particular, this extends the range treated in \cite{lcc1} to $4\le p<5$. In the
overlapping range $3<p<4$, our a priori estimate also shows that the translational regularity required in \cite{lcc1} follows automatically for weak solutions.

\subsection*{The cubic case $p=3$}
For $p=3$, the nonlinearity satisfies $|f(u)|\lesssim |u|+|u|^3$, and since $u\in H^1(\R^3)\hookrightarrow L^6(\R^3)$, we have $f(u)\in L^2(\R^3)$. Thus, we can rewrite \eqref{eq:problem} as
$$
  u_{tt}+\ga u_t-\Delta u + m u - g = -f(u)\in L^2(\R^3).
$$
The forcing term is incorporated into the potential energy; see \cite{lcc} for the corresponding argument on bounded domains. Moreover, under the cubic growth bound
$|f'(s)|\lesssim 1+|s|^2$, Sobolev embedding gives
$$
  \|f(u)-f(v)\|_{L^2(\R^3)} \lesssim \bigl(1+\|u\|_{H^1(\R^3)}^2+\|v\|_{H^1(\R^3)}^2\bigr) \|u-v\|_{H^1(\R^3)}.
$$
For two weak solutions with the same forcing, the difference equation therefore admits the standard energy estimate. Gronwall's inequality then yields unconditional uniqueness in $\H$ on $\R^3$.

\subsection*{The sub-quartic case $3<p<4$}
For $p>3$ the term $f(u)$ need no longer lie in $L^2(\R^3)$ and the reformulation above fails. The key observation is that, although $f(u)$ itself may not be square-integrable, the difference of the nonlinearities at two solutions $u,v$ is controlled by H\"older's inequality:
\begin{equation}\label{eq:u-v}
  \|f(u)-f(v)\|_{L^2(\R^3)}
  \le C\|u-v\|_{L^2(\R^3)}
  +C\l(\|u\|_{L^6(\R^3)}^{p-1}+\|v\|_{L^6(\R^3)}^{p-1}\r)
  \|u-v\|_{L^{\frac{6}{4-p}}(\R^3)},
\end{equation}
valid for $3<p<4$. Hence, if $u-v$ has the additional integrability $u-v\in L^{6/(4-p)}(\R^3)$, then $f(u)-f(v)\in L^2(\R^3)$, allowing energy estimates for the difference equation. To accommodate both these estimates and the attractor analysis, in \cite{lcc1} we introduced \emph{translationally regular solutions}:
those solutions $u$ for which the difference $\delta=u-h$ from some stationary solution $h$ satisfies
$$
  \delta\in L_t^{4/(p-2)}L_x^{12/(4-p)}(\R^3).
$$
For such solutions the nonlinear difference acts as an $L^2$-source, and an iteration argument yields global well-posedness and continuous dependence. On a bounded domain the Strichartz estimates lose derivatives \cite{sogge09}, which caps the growth at $p\le 25/7$, whereas on $\R^3$ the sharp estimates cover the whole range $3<p<4$. When $g\in L^2$ one reaches the quintic exponent $p\le 5$ because the solution gains the extra regularity $u\in L_t^4L_x^{12}$, whence
$$
  \|u^5\|_{L_t^1L_x^2}\le \|u\|_{L_t^4L_x^{12}}^4\|u\|_{L_t^\infty L_x^6};
$$
see \cite{zelik16} for details.
\subsection*{The sub-quintic regime $4\le p<5$}
For $p>4$ the exponent $6/(4-p)$ is negative, while at $p=4$ the limiting estimate requires spatial $L^\infty$ control; hence \eqref{eq:u-v} cannot be used within the preceding Strichartz framework in the range $4\le p<5$. We refine the treatment of the nonlinear difference by expanding around a stationary solution. Let $h\in H^1(\R^3)$ solve
\begin{equation}\label{eq:heq}
  -\Delta h + m h + f(h)=g,\qquad g\in H^{-1}(\R^3),
\end{equation}
whose existence follows from standard variational arguments (see Lemma \ref{le_sta_sol}). For a solution $u(t)$ of \eqref{eq:problem}, define the perturbation
$$
  \delta(t):=u(t)-h.
$$
Subtracting \eqref{eq:heq} from \eqref{eq:problem}, we obtain
\begin{equation}\label{eq:delta_intro}
  \delta_{tt}+\ga\delta_t-\Delta\delta + m\delta + f(u)-f(h)=0.
\end{equation}
We expand the super-cubic part $f_1$ around $h$ and keep $f_2$, which is globally Lipschitz by \eqref{eq:fdecomp}:
\begin{equation}\label{eq:fu-fv}
  \begin{aligned}
    f(u)-f(h) & = f_1(h+\delta)-f_1(h)+f_2(u)-f_2(h)                                                                          \\
              & = f_1^{\prime}(h)\delta + \delta^2 \int_0^1 (1-\tau) f_1^{\prime\prime}(h+\tau\delta)\,\td\tau+f_2(u)-f_2(h).
  \end{aligned}
\end{equation}
Define
$$
  V(x):=f_1^{\prime}(h(x)),\qquad
  \Gamma(t,x):= \int_0^1 (1-\tau) f_1^{\prime\prime}(h(x)+\tau \delta(t,x))\,\td\tau,
$$
so that \eqref{eq:fu-fv} becomes
\begin{equation}\label{eq:VGamma}
  f(u)-f(h) = V\delta + \delta^2\Gamma+f_2(u)-f_2(h).
\end{equation}
Substituting this into \eqref{eq:delta_intro} yields a wave equation with a potential,
\begin{equation}\label{eq:delta_V}
  \delta_{tt}+\ga\delta_t -\Delta\delta + V\delta + m\delta =-\delta^2\Gamma+f_2(h)-f_2(u),
\end{equation}
whose linear part is governed by the Schr\"odinger operator $\sop:=-\Delta+V$.

\noindent \textbf{Key observations and advantages of this decomposition.}

\begin{enumerate}
  \item \textit{The potential $V=f_1^{\prime}(h)$ belongs to the Kato class $\K_0$.} Since $f_1^{\prime\prime}(0)=0$ and $|f_1^{\prime\prime\prime}(s)|\le C|s|^{p-3}$, we have $|f_1^{\prime\prime}(s)|\le C|s|^{p-2}$, so that $\nabla V=f_1^{\prime\prime}(h)\nabla h$ obeys
        $$
          \|\nabla V\|_{L^1}\le C\|h\|_{L^{2(p-2)}(\R^3)}^{p-2}\|\nabla h\|<\infty,
        $$
        because $2(p-2)\in(2,6)$ for $3<p<5$ and $h\in H^1(\R^3)\hookrightarrow L^2(\R^3)\cap L^6(\R^3)$. Hence, $V\in \dot W^{1,1}(\R^3)\subset \K_0$ (see Subsection~\ref{sec:kato}). This is exactly why the super-cubic part $f_1$ is split off: for $f$ itself, using the bound $|f^{\prime\prime}(s)|\le C(1+|s|^{p-2})$ in this argument would require $\nabla h\in L^1(\R^3)$, which is not guaranteed by $h\in H^1(\R^3)$. In \eqref{eq:V-egy} below, with the cut-off parameter $\Lambda$ in \eqref{eq:fdecomp_def} chosen sufficiently large, the potential $V=f_1^{\prime}(h)$ itself has arbitrarily small Kato norm, and under this smallness condition $\L$ has purely absolutely continuous spectrum and no resonances on $[0,+\infty)$. This is what produces the spectral properties needed for the dispersive estimates.

  \item \textit{The linear part of \eqref{eq:delta_V} is governed by $\L$.} Equation \eqref{eq:delta_V} reads
        $$ \delta_{tt}+\ga\delta_t + \L\delta + m\delta = -\delta^2\Gamma+f_2(h)-f_2(u), $$
        so the linear evolution can be treated using the wave propagator of $\L$, for which Sections~\ref{sec:besov}--\ref{sec:strichartz} provide Strichartz estimates in the associated Besov scale.

  \item \textit{The remainder $\delta^2\Gamma$ is quadratic in $\delta$.} Its coefficient $\Gamma$ is bounded in the relevant Besov norms, since $f_1^{\prime\prime}(h+\tau\delta)$ is controlled by the $H^1$-norm of $h+\tau\delta$ and the growth of $f_1^{\prime\prime\prime}$. The product rules of Propositions~\ref{paraproduct}, \ref{negative_index} and~\ref{comp} give
        $$
          \|\delta^2\Gamma\|_{\dot B_{r',2}^{\lambda_k}(\R^3)} \lesssim R^{p-1}\|\delta\|_{\dot B_{r,2}^{\rho_{k-1}}(\R^3)},
        $$
        a gain of regularity that lets an iteration argument close. This is where the Strichartz estimates for the perturbed operator are essential to the present iteration.

  \item \textit{Why expand around a stationary solution?} Taking instead the difference of two arbitrary solutions $u$ and $v$ would produce $f_1^{\prime}(v)$ as the potential, hence a time-dependent operator, for which the spectral theory used here is unavailable. Choosing $h$ stationary is the key step in our argument for reaching the range $3< p<5$.
\end{enumerate}

A natural question is why we use Besov spaces in the present argument, rather than working only with the classical Sobolev spaces. There are two advantages. First, for super-cubic nonlinearities, the lack of 
parabolic smoothing and the possible failure of $L^2$-integrability of $f(u)$ require a precise description of the regularity of the nonlinear difference term. The iterative argument involves fractional
derivatives and negative regularity indices, and the fine scale of Besov spaces provides the product and composition estimates needed to close the argument. Second, the dyadic decomposition generated by
$\sqsop$ is directly compatible with the perturbed wave propagator. We therefore first establish the linear estimates in the spectral Besov spaces associated with $\sqsop$. The norm equivalence in
Theorem~\ref{eqv_besov} then allows these estimates to be transferred to the classical homogeneous Besov spaces. 

To handle this perturbed linear operator, we use the dispersive theory for wave equations with Kato-class potentials. For potentials in this class, $L^p$ estimates and dispersive $L^\infty$ decay estimates were proved under smallness or sign assumptions in \cite{beals93,Dancona05,georgiev,Pierfelice07}. Later works developed estimates for scaling-critical classes of potentials using a structure theory of wave operators; see \cite{beceanu14,Beceanu}. The hypotheses in these results are not identical. In particular, the dispersive estimate in \cite{bui} assumes that the potential belongs to
$\K_0$ and that $\sop$ has no eigenvalues or resonances on the nonnegative real axis; see also \cite{Rodnianski}.

For the present problem, the nonlinearity decomposition \eqref{eq:fdecomp} yields the potential $V=f_1^\prime(h)$. By Lemma \ref{lem_V}, this potential satisfies the spectral assumptions needed
to apply \cite{bui}. Under our potential assumptions, we obtain the required equivalence between the Besov spaces associated with $\sqsop$ and the classical homogeneous Besov spaces by adapting the
argument of \cite[Proposition 3.5]{Iwabuchi}. This equivalence allows us to transfer the dispersive estimates to our setting and derive the Strichartz estimates for the perturbed wave equation. These estimates
are then used in the Besov-scale iteration underlying the a priori estimate below and in proving strong time continuity of weak solutions in the energy space $\H$.

\begin{thm}[A priori estimate and time continuity]\label{regular}
  Let $g\in H^{-1}(\R^3)$ and let the nonlinearity $f$ satisfy \eqref{F1}, \eqref{F2}. Suppose that $u$ is a weak solution of equation \eqref{eq:problem} with initial data $(u_0,u_1) \in \H$. Then, $(u,u_t) \in C([0,T];\H)$ for any $T>0$, and the a priori estimate
  \begin{equation}\label{eq:regu_est}
    \| u-h\|_{L^{\frac{8}{p-1}}(t_0,t_0+1;L^{\frac{24}{5-p}}(\R^3))}\le Q\left(\sup_{t_0\le t\le t_0+1} (\|u(t)\|_{H^1(\R^3)}+\|u_t(t)\|_{L^2(\R^3)})+\|g\|_{H^{-1}(\R^3)}\right)
  \end{equation}
  holds for any $t_0\ge 0$, where $h$ is a stationary solution and $Q$ is a continuous increasing function whose dependence is restricted to the fixed structural constants of the equation, $p,m,\gamma$, and $\|g\|_{H^{-1}(\R^3)}$.
\end{thm}

\subsection*{\textbf{Unconditional uniqueness of weak solutions}}
Theorem \ref{regular} shows that \emph{every} weak solution automatically satisfies the Strichartz bound \eqref{eq:regu_est}. The uniqueness of weak solutions (Theorem \ref{thm_uniq}) is therefore
\emph{unconditional}: it holds in the class of all weak solutions with data in $\H$, no auxiliary space-time integrability being assumed. In \cite{lcc1}, the corresponding estimate was assumed through the notion of a translationally regular solution. The a priori estimate above shows that, in the present setting, every weak solution is translationally regular. Thus, in the overlapping range $3<p<4$, the additional regularity condition used in \cite{lcc1} is automatic for weak solutions. 

The natural benchmark for the admissible range of $p$ is the undamped, unforced equation
\begin{equation}\label{eq:nlw}
  u_{tt}-\Delta u+f(u)=0 \qquad\text{on } \R^d .
\end{equation}
For \eqref{eq:nlw} the energy-subcritical range is $p<p_c:=\frac{d+2}{d-2}$ when $d\ge3$, and there unconditional uniqueness
in $H^1(\R^d)\times L^2(\R^d)$ goes back to Ginibre and Velo \cite{Ginibre}, whose argument covers all $d\ge2$. The critical exponent $p=p_c$ is considerably more delicate and has been reached only in
higher dimensions: it was first addressed by Planchon \cite{planchon03}, who established unconditional uniqueness in dimensions $d=4,5,6$. For $d\ge6$ the same conclusion was obtained by Bulut et al. \cite{bulut},
who had in addition to contend with the limited smoothness of a nonlinearity of order $p_c<3$. To the best of our knowledge, in dimension three at the quintic exponent $p=5$, unconditional uniqueness in the energy space remains open. Uniqueness is known under additional assumptions, for example in the Strichartz class (see \cite{ss}) or under an energy inequality (see \cite{masmoudi, struwe99}).

The quintic endpoint is not addressed by the present argument. Together with the cubic energy argument above, our results give unconditional uniqueness in $\H$ for $3\le p<5$ under the corresponding 
growth assumptions. At $p=5$ the present iteration meets an endpoint obstruction: the intermediate Strichartz pair reaches the forbidden wave-Strichartz endpoint $(2,\infty)$ in the present three-dimensional framework, so that the ladder of Besov norms can no longer be climbed (see Remark \ref{rem_iteration}). This is an obstruction of the present argument rather than a
statement of optimality; reaching the critical exponent would require a different mechanism.

\subsection*{\textbf{Long-time dynamics and attractors}}
Once well-posedness is available, we turn to the asymptotic behavior. On an unbounded domain, the Sobolev embeddings
are not compact, so compactness cannot be inferred from regularity alone. We proceed instead by combining a decomposition with an asymptotically regular component and uniform tail estimates
of \cite{Wang1999}. Using the weighted energy identity with the spatial cut-off $\chi_n(x)=\chi(x/n)$, we show that for every bounded set of initial data and every $\eps>0$ there are $T_0$ and $n_0$ with
$$
  \int_{|x|\ge n_0} \l(|u_t|^2+|\nabla u|^2+|u|^2\r)\,\tdx < \eps,
  \qquad t\ge T_0 .
$$
Together with dissipativity and a refined decomposition of the solutions, this tail smallness gives the asymptotic compactness of the semigroup $S(t)$ in $\H=H^1(\R^3)\times L^2(\R^3)$, and the existence of a compact global attractor follows.
\begin{thm}[Global attractor]\label{thm_attr}
  Let $g\in H^{-1}(\R^3)$, and let $f\in C^3(\R)$ satisfy \eqref{F1}--\eqref{F3}. Let $\{S(t)\}_{t\geq 0}$ be the semigroup generated by solutions of problem \eqref{eq:problem} in $\H$. Then $\{S(t)\}_{t\geq 0}$ possesses a compact global attractor $\mathscr{A}$ in $\H$.
\end{thm}
The paper is organised as follows. Section~\ref{sec:besov} provides the functional framework, introducing the Kato class and the Besov spaces associated with $\sop$ alongside their equivalence with classical Besov spaces and necessary product rules. The corresponding Strichartz estimates for the wave equation with potentials are established in Section~\ref{sec:strichartz}. In Section~\ref{sec:wp}, the well-posedness of \eqref{eq:problem} is settled through an iteration argument on the Besov scale. Finally, Section~\ref{sec:attractor} concludes with the existence of the global attractor, combining dissipativity, uniform tail estimates, and a refined decomposition of the solutions.

The symbol $C$ denotes a generic positive constant, which may change from line to line; the numbered constants $C_i$, $i\in \mathbb{N}$, retain their values once introduced.
\section{Besov spaces associated with the Schr\"odinger operator}\label{sec:besov}
\subsection{The Kato class}\label{sec:kato}
A central notion in the theory of Schr\"odinger operators is the Kato class, introduced by Kato in his work \cite{kato} on the self-adjointness of quantum Hamiltonians. A potential $V$ belongs to the Kato class $\K$ if it satisfies the \emph{Kato condition}
\begin{equation}
  \lim_{r \to 0} \, \sup_{x \in \R^n} \int_{\lvert x-y \rvert \leq r} \frac{\lvert V(y) \rvert}{\lvert x-y \rvert^{n-2}}  \td y = 0 \quad (\text{for } n \geq 3),
\end{equation}
with analogous conditions in dimensions $n=1,2$. The Kato class contains potentials outside every $L^p$-space, in particular Coulomb-type singularities $|x|^{-a}$ with $0<a<2$. It is the natural class for heat-kernel estimates and Feynman--Kac representations. The Kato class provides a useful framework for defining Schr\"odinger operators and for deriving the estimates used below; the specific spectral properties required here will be imposed and verified separately; see \cite{simon82}.

For dispersive estimates on the linear propagator of $\sop = -\Delta + V$, the appropriate class is the \emph{Kato space} $\K_0$, the closure of the bounded compactly supported functions in the global Kato norm (see \cite{Rodnianski})
$$\|V\|_{\K}=\sup_{x \in \R^n} \int_{\R^n} \frac{\lvert V(y) \rvert}{\lvert x-y \rvert^{n-2}}  \td y \quad (\text{for } n \geq 3).$$
In contrast to the small-scale condition in the definition of $\K$ above, the global Kato norm controls the weighted integral of the potential over all of $\R^n$, uniformly with respect to the center
$x$. This is the quantity that controls the dispersive and Strichartz estimates below.

The following lemma records the role of smallness of the Kato norm in obtaining the spectral properties of the perturbed operator that will be used below.
\begin{lemma}\label{lem_V}
The following properties hold:
    \begin{itemize}
    \item[(i)] Let $h\in H^1(\R^3)$ and let $f=f_1+f_2$ be as in \eqref{eq:fdecomp}. Set $V:=f_1^{\prime}(h)$. Then $V\in \dot W^{1,1}(\R^3)\cap L^{3/2}(\R^3)\subset \K_0\cap L^{3/2}(\R^3)$, and
      $$
        \|V\|_{\K}\le\|\nabla V\|_{L^1(\R^3)}\le C\|h\|^{p-1}_{H^1(\R^3)}.
      $$
    \item[(ii)] If $\|V\|_{\K}=\kappa<4\pi$, then $\sop:=-\Delta+V$ has purely absolutely continuous spectrum
      and no eigenvalues or resonances on $[0,+\infty)$. Here, a resonance at $\lambda\ge0$ means a solution $\psi$ of
      $$
        \psi+(-\Delta-\lambda\pm i0)^{-1}(V\psi)=0
      $$
      for either choice of sign, such that
      $$
        \psi\notin L^2(\R^3),\qquad (1+|x|^2)^{-s/2}\psi\in L^2(\R^3) \quad\text{for every }s>\tfrac12.
      $$
      Moreover, $P_c=\mathrm{Id}$ and
      \begin{equation}\label{eq:V_coercive}
        \int_{\R^3}\bigl(|\nabla\varphi|^2+V\varphi^2\bigr)\,\tdx
        \ge\left(1-\frac{\|V\|_{\K}}{4\pi}\right)
        \|\nabla\varphi\|^2,\qquad \varphi\in H^1(\R^3),
      \end{equation}
      whence $\sop\ge0$.
  \end{itemize}
\end{lemma}

\begin{proof}
  (i) From $f_1^\prime(0)=f_1^{\prime\prime}(0)=0$ and $|f_1^{\prime\prime\prime}(s)|\le C|s|^{p-3}$ we get
  $$|f_1^{\prime\prime}(s)|\le C|s|^{p-2},\qquad |f_1^\prime(s)|\le C|s|^{p-1}.$$
  Since $\nabla V=f_1^{\prime\prime}(h)\nabla h$, H\"older's inequality yields
  $$
    \|\nabla V\|_{L^1}\le \|f_1^{\prime\prime}(h)\|\|\nabla h\|
    \le C\|h\|^{p-2}_{L^{2(p-2)}}\|\nabla h\|,
  $$
  which is finite because $2(p-2)\in(2,6)$ for $3<p<5$ and $H^1(\R^3)\hookrightarrow L^q(\R^3)$ for
  $2\le q\le 6$.

  Similarly, since $3(p-1)/2\in(3,6)$, $$\|V\|_{L^{3/2}}\le C\|h\|^{p-1}_{L^{3(p-1)/2}}<\infty.$$
  The inclusion $\dot W^{1,1}\subset\K_0$ and the bound $\|V\|_\K\le\|\nabla V\|_{L^1}$ complete the proof of (i).

  (ii) By Kato's inequality \cite{kato66} (see also \cite{simon82}),
  $$\left|\int_{\R^3}V\varphi^2\td x\right| \le\frac{\|V\|_{\K}}{4\pi}\|\nabla\varphi\|^2,$$
  which gives \eqref{eq:V_coercive} because $\|V\|_{\K}<4\pi$.
  In particular $\sop\ge0$.

  Since $V\in\K_0$ by (i) and $\|V\|_{\K}<4\pi$, the small-Kato criterion stated in \cite[Section~1.2]{bui}
  rules out eigenvalues and resonances on $[0,+\infty)$ in the sense used there. By Kato's theorem \cite[Theorem 6.1]{kato66},
  $\sop$ is unitarily equivalent to $-\Delta$ under the same smallness condition. Thus its spectrum is purely absolutely
  continuous and $P_c=\mathrm{Id}$. Moreover, $|V|\le C|h|^{p-1}$ and $2<p-1<4$, so $h\in H^1(\R^3)$ implies $V\in L^1(\R^3)$.

\end{proof}

\subsection{Besov spaces associated with the Schr\"odinger operator}
Besov spaces provide a finely tuned measure of regularity. Although the classical Fourier-based spaces remain useful, the dyadic decomposition associated with
$\sop=-\Delta+V$, $V\in\K_0$, is more directly compatible with the operator governing the evolution. We therefore first work with the Besov spaces associated with $\sop$, and later use the norm equivalence
in Theorem~\ref{eqv_besov} to return to the classical homogeneous Besov scale.

Let $\phi \in C_c^\infty(\R)$ be non-negative, with
\begin{equation}
  \operatorname{supp} \phi=\{\lm\in \R: 2^{-1}\le \lm\le 2\}
\end{equation}
and
\begin{equation}
  \sum_{j=-\infty}^{\infty}\phi(2^{-j}\lm)=1 \quad\text{for } \lm > 0 ,
\end{equation}
and set $\phi_j(\lm)=\phi(2^{-j}\lm)$ for $j\in \mathbb{Z}$.

Define
$$
  P_c:=I-\sum_{\lambda\in\sigma_{\mathrm p}(\sop)} E_{\sop}(\{\lambda\}),
$$
where $E_{\sop}$ is the projection-valued spectral measure of the self-adjoint operator $\sop$ on $L^2(\R^3)$ and the sum converges in the strong operator topology. Thus, $P_c$ projects onto the continuous spectral subspace; in particular, $P_c=I$ if $\sop$ has no eigenvalues.

For $V\in\K_0$, the continuous spectrum of $\sop$ lies in $[0,\infty)$. Following \cite{bui}, we interpret spectral multipliers involving $\sqsop$ on $P_cL^2(\R^3)$, extending them by zero
on the orthogonal complement; in particular,
$$
  \phi_j(\sqsop)=\int_{[0,\infty)}\phi_j(\sqrt{\lambda})\,dE_{\sop}(\lambda)P_c.
$$
The same convention applies to powers of $\sop$ composed with spectral cutoffs and to the wave propagators on $P_cL^2$; the quoted $L^p$ bounds refer to their bounded extensions.

For $1\le p<\infty$, $1\le q<\infty$, and
$$
  s\in I_p:=
  \left(-\min\{2,3(1-\frac{1}{p})\},\,\min\{2,\frac{3}{p}\}\right),
$$
let
$$
  \mathcal D_{s,p,q}(\sop)
  :=\left\{f\in P_cC_c^\infty(\mathbb R^3):
  \|f\|_{\dot B^s_{p,q}(\sqrt{\sop})}<\infty\right\},
$$
where
\begin{equation}\label{eq:besov_norm}
  \|f\|_{\dot B^s_{p,q}(\sqrt{\sop})}
  :=\left(\sum_{j\in\mathbb Z}
  2^{jsq}\|\phi_j(\sqrt{\sop})f\|_{L^p(\mathbb R^3)}^q
  \right)^{1/q}.
\end{equation}
We define $\dot B^s_{p,q}(\sqrt{\sop})$ as the completion of $\mathcal D_{s,p,q}(\sop)$ in this norm. The definition is used below only in the indicated parameter range.

\subsection{Equivalence of Besov spaces}

In this subsection we state the equivalence
$$\dot{B}^s_{p,q}(\sqsop)\cong \dot{B}^s_{p,q}(\R^3)$$
between the perturbed Besov spaces built on $\sop$ and the standard ones on $\R^3$.

The equivalence rests on estimates for the spectral multipliers of $\sop$, the basic tool being the Bernstein inequality adapted to $\sop$. The $L^p$-boundedness of $\phi_j(\sqsop)$ was proved in \cite{Dancona05} under a smallness assumption on the negative part $V_-=\min(V,0)$; this restriction was later removed in \cite{bui} by means of an operator-valued Wiener theorem, which handles potentials of large Kato norm.
\begin{prop}[Bernstein inequality \cite{bui}]\label{prop_Bernstein}
  Suppose that the potential $V\in \K_0$ and $\sop$ has no eigenvalues or resonances on the nonnegative real axis $[0, +\infty)$. Let $\varphi \in C^{\infty}_c([0, +\infty))$ be a compactly supported smooth function whose support is contained in the interval $[0, 2]$. Then, if $1 \le p\le q\le \infty$, there exists a constant $C>0$ such that
  \begin{equation}\label{eq:Bernstein}
    \|\varphi(2^{-j}\sqsop)f\|_{L^{q}(\R^3)}\le C2^{3(\frac1{p}-\frac1q)j}\|f\|_{L^{p}(\R^3)},
  \end{equation}
  for any $j\in \Z$ and $f\in L^{p}(\R^3)$.
\end{prop}

\begin{cor}\label{cor_bern}
  Let the potential $V$ and the function $\varphi$ be as in Proposition \ref{prop_Bernstein}. Then, if $\al\ge 0$ and $1 \le p\le q\le \infty$, there exists a constant $C>0$ such that
  \begin{equation}\label{eq:Hbern}
    \|\sop^{\al}\varphi(2^{-j}\sqsop)f\|_{L^{q}(\R^3)}\le C2^{3(\frac1{p}-\frac1q)j+2\al j}\|f\|_{L^{p}(\R^3)},
  \end{equation}
  for any $j\in \Z$ and $f\in L^{p}(\R^3)$.
\end{cor}

Suppose that $V\in\K_0\cap L^{3/2}(\R^3)$, $\|V_-\|_{\K}<\pi$, and $\sop$ has no eigenvalues or resonances on $[0,\infty)$. The embedding $L^{3/2}(\R^3)\hookrightarrow L^{3/2,\infty}(\R^3)$ and the Bernstein inequality \eqref{eq:Bernstein} ensure that the hypotheses required in the proof of \cite[Proposition 3.5]{Iwabuchi} are satisfied. Adapting that argument to the present setting yields the following norm equivalence.
\begin{thm}\label{eqv_besov}
  Let $V\in\K_0\cap L^{3/2}(\mathbb R^3)$ satisfy $\|V_-\|_{\K}<\pi$, and assume that $\sop$ has no eigenvalues or
  resonances on $[0,+\infty)$. Let $1\le p<\infty$, $1\le q<\infty$, and $s\in I_p$. Then
  $$
    \dot B^s_{p,q}(\sqrt{\sop}) \cong \dot B^s_{p,q}(\mathbb R^3),
  $$
  with equivalent norms.
\end{thm}
\begin{proof}
Under the present assumptions, $\sop\ge0$ and $P_c=I$.
The norm comparison described above gives
$$
  \|f\|_{\dot B^s_{p,q}(\sqrt{\sop})} \cong \|f\|_{\dot B^s_{p,q}(\mathbb R^3)},\qquad f\in C_c^\infty(\mathbb R^3).
$$
Since $p,q<\infty$ and $-3(1-1/p)<s<3/p$, this common core is dense in the classical homogeneous Besov space. Hence the identity on the core extends uniquely to an isomorphism between the corresponding completions. We use the canonical distributional realization of the classical space: as $s<3/p$, the free Littlewood--Paley reconstruction converges in $\mathcal S'(\mathbb R^3)$ and selects the representative whose low-frequency cutoffs tend to zero in $\mathcal S'(\mathbb R^3)$. Thus the completion used here has a unique distributional realization, with no polynomial ambiguity.
\end{proof}

\begin{remark}\label{constant_smallV}
Let $\mathcal V\subset\K_0\cap L^{3/2}(\R^3)$ be a family of real-valued
potentials such that
$$
  \sup_{V\in\mathcal V}\|V\|_{\K}\le\kappa<\pi,
  \qquad
  \sup_{V\in\mathcal V}\|V\|_{L^{3/2}(\R^3)}\le M,
$$
and assume that $-\Delta+V$ has no eigenvalues or resonances on $[0,+\infty)$ for every $V\in\mathcal V$. For the Bernstein estimate, the small-potential proof of \cite[Proposition 2.2 and Lemma 2.4]{bui} gives, for fixed $1\le a\le b\le\infty$ and sufficiently large $m$, the constant in the Bernstein inequality \eqref{eq:Bernstein} is bounded by
$$C_{a,b,m}\|\varphi\|_{W^{m,1}([0,+\infty))}\left(1-\frac{\|V\|_{\K}}{4\pi}\right)^{-2}.$$
Thus, for fixed $\varphi$, this estimate is uniform over $V\in\mathcal V$.

Moreover, the $L^{3/2}$ bound gives a uniform $L^{3/2,\infty}$ bound, while $\|V\|_{\K}\le\kappa<\pi$ gives the uniform Kato gap
$\pi-\|V_-\|_{\K}\ge\pi-\kappa>0$. Hence, for fixed admissible indices $r,q,s$, the proof of \cite[Proposition 3.5, Section 9]{Iwabuchi}
yields
$$
  C^{-1}\|\varphi\|_{\dot B^s_{r,q}(\R^3)}
  \le
  \|\varphi\|_{\dot B^s_{r,q}(\sqsop)}
  \le
  C\|\varphi\|_{\dot B^s_{r,q}(\R^3)}
  \qquad (V\in\mathcal V),
$$
where $C$ is independent of $V$ and may be chosen to depend only on the fixed indices, the structural constants, $M$, and $\pi-\kappa$.
\end{remark}

\subsection{Product rules}

Nonlinear problems require estimates for products and compositions in Besov spaces. We begin with the product rules.

We first record the boundedness of products in homogeneous Besov spaces of equal regularity.

\begin{prop}\label{paraproduct}
  Let $s>0$ and let
  $$
    1\le p_1,p_2<\frac{n}{s},\qquad
    \frac{1}{p_0}
    =
    \frac{1}{p_1}+\frac{1}{p_2}-\frac{s}{n},
    \qquad
    1\le p_0<\infty .
  $$
  For $i=1,2$, define $p_i^*$ by
  $$
    \frac{1}{p_i^*}=\frac{1}{p_i}-\frac{s}{n},
  $$
  and assume that $p_i^*\ge 2$. Then
  $$
    \|fg\|_{\dot B^s_{p_0,2}(\R^n)}
    \lesssim
    \|f\|_{\dot B^s_{p_1,2}(\R^n)}
    \|g\|_{\dot B^s_{p_2,2}(\R^n)} .
  $$
\end{prop}

\begin{proof}
  By the homogeneous Sobolev embedding,
  $$
    \dot B^s_{p_i,2}(\R^n)
    \hookrightarrow L^{p_i^*}(\R^n),
    \qquad i=1,2,
  $$
  where the assumption $p_i^*\ge2$ ensures that the fine index $2$ is admissible. Hence,
  $$
    \|f\|_{L^{p_1^*}(\R^n)}\lesssim\|f\|_{\dot B^s_{p_1,2}(\R^n)},\qquad\|g\|_{L^{p_2^*}}\lesssim\|g\|_{\dot B^s_{p_2,2}(\R^n)}.
  $$

  Let
  $$
    fg=\dot T_f g+\dot T_g f+\dot R(f,g)
  $$
  be the homogeneous Bony decomposition. Since
  $$
    \frac{1}{p_0}=\frac{1}{p_1}+\frac{1}{p_2^*}=\frac{1}{p_1^*}+\frac{1}{p_2},
  $$
  the homogeneous paraproduct and remainder estimates of
  \cite[Theorems 2.47 and 2.52]{bcd} give
  $$
    \begin{aligned}
      \|\dot T_g f\|_{\dot B^s_{p_0,2}(\R^n)}+\|\dot R(f,g)\|_{\dot B^s_{p_0,2}(\R^n)}
                                              & \lesssim\|g\|_{L^{p_2^*}(\R^n)}\|f\|_{\dot B^s_{p_1,2}(\R^n)}, \\
      \|\dot T_f g\|_{\dot B^s_{p_0,2}(\R^n)} & \lesssim\|f\|_{L^{p_1^*}}\|g\|_{\dot B^s_{p_2,2}(\R^n)} .
    \end{aligned}
  $$
  Therefore,
  $$
    \begin{aligned}
      \|fg\|_{\dot B^s_{p_0,2}(\R^n)}
       & \lesssim
      \|f\|_{\dot B^s_{p_1,2}(\R^n)}\|g\|_{L^{p_2^*}(\R^n)}+\|f\|_{L^{p_1^*}(\R^n)}\|g\|_{\dot B^s_{p_2,2}(\R^n)} \\
       & \lesssim\|f\|_{\dot B^s_{p_1,2}(\R^n)}\|g\|_{\dot B^s_{p_2,2}(\R^n)} .
    \end{aligned}
  $$
  This proves the result.
\end{proof}
The case of mixed regularity, one factor having negative smoothness, is particularly important for us: a negative-order Besov norm still controls the product, at the cost of regularity in the target space.

\begin{prop}\label{negative_index}
  Let $s_1,s_2\in\R$ and $1\leq p_0,p_1,p_2<\infty$ satisfy
  $$
    \frac{1}{p_0}=\frac{1}{p_1}+\frac{1}{p_2},\qquad 0<s_1<\frac{n}{p_1},\qquad s_2<0,\qquad s_1+s_2>0.
  $$
  Let $f$ be a tempered distribution, represented in a fixed realization of the homogeneous Besov spaces, such that
  $$
    f\in \dot B^{s_1}_{p_1,\infty}(\R^n)\cap \dot B^{s_2}_{p_2,2}(\R^n).
  $$
  Then the product $f^2$, defined by the corresponding Littlewood--Paley paraproduct decomposition, belongs to
  $\dot B^{s_1+s_2}_{p_0,2}(\R^n)$ and satisfies
  \begin{equation}\label{eq:negative_index}
    \|f^2\|_{\dot B^{s_1+s_2}_{p_0,2}(\R^n)} \lesssim \|f\|_{\dot B^{s_1}_{p_1,\infty}(\R^n)} \|f\|_{\dot B^{s_2}_{p_2,2}(\R^n)}.
  \end{equation}
  The implicit constant depends only on the dimension, the indices, and the choice of the Littlewood--Paley decomposition.
\end{prop}

\begin{proof}
  Fix a homogeneous Littlewood--Paley decomposition $\{\Delta_j\}_{j\in\mathbb Z}$ and write
  $$
    S_{j-1}f:=\sum_{k\leq j-2}\Delta_k f, \qquad \widetilde\Delta_j f:=\sum_{|k-j|\leq 1}\Delta_k f.
  $$
  Bony's decomposition gives
  \begin{equation}\label{eq:bony-square}
    f^2=2\sum_{j\in\mathbb Z}(S_{j-1}f)\Delta_j f+\sum_{j\in\mathbb Z}(\Delta_j f)(\widetilde\Delta_j f).
  \end{equation}
  The first sum contains the low--high interactions, while the second contains the interactions between comparable frequencies.

  Set
  $$
    a_j:=2^{s_1j}\|\Delta_j f\|_{L^{p_1}(\R^n)},\qquad b_j:=2^{s_2j}\|\Delta_j f\|_{L^{p_2}(\R^n)}.
  $$
  Then
  $$
    \|a\|_{\ell^\infty}=\|f\|_{\dot B^{s_1}_{p_1,\infty}(\R^n)},\qquad\|b\|_{\ell^2}=\|f\|_{\dot B^{s_2}_{p_2,2}(\R^n)}.
  $$
  We estimate the two sums in \eqref{eq:bony-square} separately.

  For the low--high term, the Fourier support of $(S_{j-1}f)\Delta_j f$ is contained in an annulus of size $2^j$. Thus, for some fixed integer $C$,
  $$
    \Delta_\ell\bigl((S_{j-1}f)\Delta_j f\bigr)=0 \quad\text{unless}\quad |\ell-j|\leq C.
  $$
  By H\"older's inequality and the relation $1/p_0=1/p_1+1/p_2$,
  $$
    \begin{aligned}
      \|(S_{j-1}f)\Delta_j f\|_{L^{p_0}(\R^n)} & \leq \|S_{j-1}f\|_{L^{p_2}(\R^n)} \|\Delta_j f\|_{L^{p_1}(\R^n)}                               \\
                                               & \leq \left(\sum_{k\leq j-2}\|\Delta_k f\|_{L^{p_2}(\R^n)}\right)\|\Delta_j f\|_{L^{p_1}(\R^n)} \\
                                               & = 2^{-s_1j}a_j \sum_{k\leq j-2}2^{-s_2k}b_k.
    \end{aligned}
  $$
  Since $|\ell-j|\leq C$, multiplication by $2^{(s_1+s_2)\ell}$ and summation over the finitely many relevant $j$ yield
  $$
    2^{(s_1+s_2)\ell} \left\|\Delta_\ell\l(\sum_j(S_{j-1}f)\Delta_j f\r)\right\|_{L^{p_0}(\R^n)}\lesssim \|a\|_{\ell^\infty} \sum_{k\leq \ell+C} 2^{s_2(\ell-k)}b_k.
  $$
  Because $s_2<0$, the sequence
  $$
    2^{s_2m}\mathbf 1_{\{m\geq -C\}}
  $$
  belongs to $\ell^1(\mathbb Z)$. Young's inequality for sequences therefore gives
  \begin{equation}\label{eq:low-high-bound}
    \left\| \sum_j(S_{j-1}f)\Delta_j f \right\|_{\dot B^{s_1+s_2}_{p_0,2}(\R^n)} \lesssim \|a\|_{\ell^\infty}\|b\|_{\ell^2}.
  \end{equation}

  For the comparable-frequency term, the Fourier support of $(\Delta_j f)(\widetilde\Delta_j f)$ is contained in a ball of size $2^j$. Consequently,
  $$
    \Delta_\ell\bigl((\Delta_j f)(\widetilde\Delta_j f)\bigr)=0 \quad\text{if}\quad j<\ell-C.
  $$
  Again using H\"older's inequality, and noting that $\widetilde\Delta_j f$ is a sum of a fixed number of neighboring dyadic pieces, we obtain
  $$
    \left\| \Delta_\ell \sum_j(\Delta_j f)(\widetilde\Delta_j f) \right\|_{L^{p_0}(\R^n)} \lesssim \sum_{j\geq \ell-C} \|\Delta_j f\|_{L^{p_1}(\R^n)} \|\widetilde\Delta_j f\|_{L^{p_2}(\R^n)}.
  $$
  Hence,
  $$
    2^{(s_1+s_2)\ell}\left\|\Delta_\ell\sum_j(\Delta_j f)(\widetilde\Delta_j f)\right\|_{L^{p_0}(\R^n)}\lesssim\sum_{j\geq \ell-C}2^{-(s_1+s_2)(j-\ell)}a_j\sum_{|k-j|\leq 1}b_k.
  $$
  Since $s_1+s_2>0$, the sequence
  $$
    2^{-(s_1+s_2)m}\mathbf 1_{\{m\geq -C\}}
  $$
  belongs to $\ell^1(\mathbb Z)$. Also,
  $$
    \left\|\left\{a_j\sum_{|k-j|\leq 1}b_k\right\}_{j\in\mathbb Z}\right\|_{\ell^2}\lesssim\|a\|_{\ell^\infty}\|b\|_{\ell^2}.
  $$
  A second application of Young's inequality thus gives
  \begin{equation}\label{eq:high-high-bound}
    \left\|\sum_j(\Delta_j f)(\widetilde\Delta_j f)\right\|_{\dot B^{s_1+s_2}_{p_0,2}(\R^n)}\lesssim\|a\|_{\ell^\infty}\|b\|_{\ell^2}.
  \end{equation}

  The two bounds \eqref{eq:low-high-bound} and \eqref{eq:high-high-bound}, together with \eqref{eq:bony-square}, prove \eqref{eq:negative_index}.
  The estimates also ensure convergence of the paraproduct sums in the target Besov space, so the product is well-defined in the stated realization.
\end{proof}
Besides products, we need to control compositions, which is the content of the next proposition.
\begin{prop}[{\cite[Lemma 3.2]{Ginibre}}]\label{comp}
  Let $f\in C^1(\R)$ with $f(0)=0$ and $|f^{\prime} (s)|\le C|s|^{a}$ for some $a > 0$. Let $0<s<1$, let $1\le l\le k\le \infty$, $1\le q\le \infty$ and $1/\lambda=1/l-1/k$. Then the following inequality holds
  \begin{equation}
    \|f(\varphi)\|_{\dot{B}^{s}_{l,q}(\R^3)} \lesssim \|~|\varphi|^a\|_{L^{\lambda}(\R^3)}\|\varphi\|_{\dot{B}^{s}_{k,q}(\R^3)},
  \end{equation}
  for all $\varphi$ such that the norms on the right-hand side are finite.
\end{prop}

\section{Strichartz estimates for the wave equation with a potential}\label{sec:strichartz}
For the wave equation with a potential
\begin{equation}
  u_{tt} +\sop u = F, \quad u(0) = u_0, \quad u_t(0) = u_1,
\end{equation}
the continuous-spectrum part of the solution is given by the formula
\begin{equation}\label{eq:rep}
  P_cu(t) = \cos(t\sqsop)P_c u_0 + \frac{\sin(t\sqsop)}{\sqsop}P_c u_1 + \int_0^t \frac{\sin((t-\tau)\sqsop)}{\sqsop} P_c F(\tau) d\tau.
\end{equation}

The derivation of Strichartz estimates relies fundamentally on the dispersive estimate
\begin{equation}\label{eq:dispersive_estimate}
  \left\|\frac{e^{\mathrm{i}t\sqsop}}{\sqsop}P_c\varphi\right\|_{L^\infty(\R^3)} \lesssim |t|^{-1} \|P_c\varphi\|_{\dot{B}^{1}_{1,1}(\sqsop)}.
\end{equation}
For $V=0$ this is classical; the presence of a potential is what makes \eqref{eq:dispersive_estimate} delicate.

Early decay estimates for wave equations with potentials were obtained by Beals and Strauss \cite{beals93} in dimensions $n\ge 3$. Subsequent work developed Kato-norm methods for rough potentials. D'Ancona and Pierfelice \cite{Dancona05} considered potentials $V=V_1+V_2$, where $V_1$ is compactly supported with finite Kato norm, $V_2$ has small Kato norm, and $\|V_-\|_{\K}<2\pi$. Pierfelice \cite{Pierfelice07} later established a related estimate for $\sin(t\sqsop)/\sqsop$, with data in the classical inhomogeneous Besov space $B^1_{1,1}(\R^3)$, under $\|V\|_{\K}<4\pi$. Bui, Duong and Hong \cite{bui} removed the Kato-norm smallness assumptions for $V\in\K_0$, provided that $\sop$ has no eigenvalues or resonances on $[0,\infty)$. Their frequency-localised half-wave estimate \cite[Lemma 4.1]{bui}, together with the dyadic summation in the proof of their Theorem~1.1, yields \eqref{eq:dispersive_estimate} in the spectral homogeneous Besov scale.

Following Ginibre and Velo \cite{ginibre95}, who combine dispersive estimates with energy inequalities through a $TT^*$ argument, we now prove the Strichartz estimates for the wave equation with a potential.
\begin{thm}\label{thm_strichartz}
  Suppose that the potential $V\in \K_0$ and $\sop$ has no eigenvalues or resonances on the nonnegative real axis $[0, +\infty)$. Let $q_i\in(2,\infty]$, $r_i\in[2,\infty)$, $i=1,2$, satisfy
  \begin{equation*}
    -\frac12<\sigma<\frac32,
    \qquad
    -\min\left\{2,3\left(1-\frac1{r_i}\right)\right\}
    <\rho_i<\frac3{r_i},
  \end{equation*}
  and let the following admissible conditions be satisfied
  \begin{equation}\label{eq:wave_admissible}
    \begin{aligned}
       & \frac1{q_i}\le \frac12-\frac1{r_i},\quad q_i>2,                                                 \\
       & \rho_1+3(\frac12-\frac1{r_1})-\frac1{q_1}=\sigma=1-(\rho_2+3(\frac12-\frac1{r_2})-\frac1{q_2}).
    \end{aligned}
  \end{equation}
  Then the following Strichartz estimate holds
  \begin{equation}\label{eq:str3}
    \|P_cu\|_{L^{q_1}(0,T;\dot B^{\rho_1}_{r_1,2}(\sqsop))}\lesssim \|P_cu_0\|_{\dot B^\sigma_{2,2}(\sqsop)}+\|P_cu_1\|_{\dot B^{\sigma-1}_{2,2}(\sqsop)}+\|P_cF\|_{L^{q_2'}(0,T;\dot B^{-\rho_2}_{r_2',2}(\sqsop))}.
  \end{equation}
\end{thm}
\begin{proof}
  By the solution representation \eqref{eq:rep}, it suffices to establish the homogeneous estimate
  \begin{equation}\label{eq:str4}
    \left\|e^{\mathrm{i}t\sqsop}P_cu_0
    \right\|_{L^{q_1}(0,T;
    \dot{B}^{\rho_1}_{r_1,2}(\sqsop))}
    \lesssim
    \|P_cu_0\|_{\dot{B}^{\sigma}_{2,2}(\sqsop)}
  \end{equation}
  and the inhomogeneous estimate
  \begin{equation}\label{eq:str5}
    \left\|\int_0^t
    \frac{\sin((t-\tau)\sqsop)}{\sqsop}
    P_cF(\tau)\,\td\tau
    \right\|_{L^{q_1}(0,T;
    \dot{B}^{\rho_1}_{r_1,2}(\sqsop))}
    \lesssim
    \|P_cF\|_{L^{q_2'}(0,T;
    \dot{B}^{-\rho_2}_{r_2',2}(\sqsop))}.
  \end{equation}

By Sobolev embedding, it suffices to prove \eqref{eq:str4} and \eqref{eq:str5} for the boundary admissible pairs
  $1/r_i=1/2-1/q_i$, $i=1,2$, for which $\rho_1=\sigma-2/q_1$ and $\rho_2=1-\sigma-2/q_2$. The remaining admissible pairs follow by Sobolev embedding, with the spatial exponents adjusted according to \eqref{eq:wave_admissible}.

  We first derive a frequency-localised dispersive estimate.
  Almost orthogonality and Corollary \ref{cor_bern} give
  $$
    \begin{aligned}
      \left\|\sqsop\phi_j(\sqsop)P_cF(\tau) \right\|_{\dot{B}^{1}_{1,1}(\sqsop)}&=\sum_{|k-j|\leq1}2^k\left\|\phi_k(\sqsop)\sqsop\phi_j(\sqsop)P_cF(\tau)\right\|_{L^1(\R^3)}\\
      &\lesssim\sum_{|k-j|\leq1}2^{2k}\|\phi_j(\sqsop)P_cF(\tau)\|_{L^1(\R^3)} \\
      &\lesssim 2^{2j}\|\phi_j(\sqsop)P_cF(\tau)\|_{L^1(\R^3)}.
    \end{aligned}
  $$
  Applying \eqref{eq:dispersive_estimate} to
  $\sqsop\phi_j(\sqsop)P_cF(\tau)$ therefore yields
  \begin{equation}\label{eq:local_dispersion}
    \left\|e^{\mathrm{i}(t-\tau)\sqsop}\phi_j(\sqsop)P_cF(\tau)\right\|_{L^\infty(\R^3)}\lesssim|t-\tau|^{-1}2^{2j}\|\phi_j(\sqsop)P_cF(\tau)\|_{L^1(\R^3)}.
  \end{equation}
  Using a spectral cutoff equal to one on the support of $\phi_j$, we interpolate this estimate with the $L^2$ unitarity of the propagator on the continuous spectral
  subspace. Since $1/r_1=1/2-1/q_1$, this gives
  \begin{equation}\label{eq:r_decay}
    \left\|e^{\mathrm{i}(t-\tau)\sqsop}\phi_j(\sqsop)P_cF(\tau)\right\|_{L^{r_1}(\R^3)}\lesssim |t-\tau|^{-2/q_1}2^{4j/q_1} \|\phi_j(\sqsop)P_cF(\tau)\|_{L^{r_1'}(\R^3)}.
  \end{equation}
  Multiplying by $2^{-2j/q_1}$ and taking the $\ell^2_j$ norm, we obtain
  $$
    \left\|e^{\mathrm{i}(t-\tau)\sqsop}P_cF(\tau)
    \right\|_{\dot{B}^{-2/q_1}_{r_1,2}(\sqsop)}
    \lesssim
    |t-\tau|^{-2/q_1}
    \|P_cF(\tau)\|_{\dot{B}^{2/q_1}_{r_1',2}(\sqsop)}.
  $$

  For $2<q_1<\infty$, Minkowski's integral inequality and the
  Hardy--Littlewood--Sobolev inequality in time, applied after
  extension of the source by zero, imply
  \begin{equation}\label{eq:rho11}
    \left\|\int_0^t e^{\mathrm{i}(t-\tau)\sqsop}
    P_cF(\tau)\,\td\tau
    \right\|_{L^{q_1}(0,T;
    \dot{B}^{\rho_1-\sigma}_{r_1,2}(\sqsop))}
    \lesssim
    \|P_cF\|_{L^{q_1'}(0,T;
    \dot{B}^{\sigma-\rho_1}_{r_1',2}(\sqsop))}.
  \end{equation}
  When $q_1=\infty$, the boundary relations give
  $r_1=2$ and $\rho_1=\sigma$. Unitarity and Minkowski's
  integral inequality instead give
  \begin{equation}\label{eq:rho11_inf}
    \left\|\int_0^t e^{\mathrm{i}(t-\tau)\sqsop}
    P_cF(\tau)\,\td\tau
    \right\|_{L^\infty(0,T;
    \dot{B}^{0}_{2,2}(\sqsop))}
    \lesssim
    \|P_cF\|_{L^1(0,T;
    \dot{B}^{0}_{2,2}(\sqsop))}.
  \end{equation}

  We next obtain the homogeneous estimate by duality. The accompanying energy estimate is carried out entirely in the Hilbert space $\dot B^0_{2,2}(\sqsop)$, whose inner product is
  \begin{equation}\label{eq:B0_inner}
    (f,g)_{\dot B^0_{2,2}(\sqsop)}:=\sum_{j\in\mathbb Z}\l(\phi_j(\sqsop)P_cf,\phi_j(\sqsop)P_cg\r)_{L^2(\R^3)} .
  \end{equation}
  Since $e^{\mathrm{i}t\sqsop}$ commutes with $\phi_j(\sqsop)$ and is unitary on $L^2(\R^3)$, it is also unitary on
  $\dot B^0_{2,2}(\sqsop)$:
  \begin{equation}\label{eq:B0_unitary}
    \l(e^{\mathrm{i}t\sqsop}P_c f,e^{\mathrm{i}t\sqsop}P_c g\r)_{\dot B^0_{2,2}(\sqsop)}=(f,g)_{\dot B^0_{2,2}(\sqsop)} .
  \end{equation}
  Using \eqref{eq:B0_unitary},
  \begin{equation}\label{eq:energy_B0}
    \begin{aligned}
        & \left\|\int_{0}^T e^{-\mathrm{i}t\sqsop}P_c F(t)\td t\right\|_{\dot B^0_{2,2}(\sqsop)}^2            \\
      = & \l(\int_{0}^T e^{-\mathrm{i}\tau\sqsop}P_c F(\tau)\td \tau,\int_{0}^T e^{-\mathrm{i}t\sqsop}P_c F(t)\td t\r)
      _{\dot B^0_{2,2}(\sqsop)}                                                                               \\
      = & \int_{0}^T
      \l(e^{\mathrm{i}t\sqsop}P_c\int_{0}^T e^{-\mathrm{i}\tau \sqsop}P_c F(\tau)\td \tau,P_c F(t)\r)
      _{\dot B^0_{2,2}(\sqsop)}\td t                                                                          \\
      = & 2\operatorname{Re}\int_{0}^{T}\int_0^t
      \l(e^{\mathrm{i}(t-\tau)\sqsop}P_c F(\tau),P_c F(t)\r)
      _{\dot B^0_{2,2}(\sqsop)}
      \td \tau\td t .
    \end{aligned}
  \end{equation}
  By H\"older's inequality on each dyadic block and Cauchy--Schwarz in the summation over $j$, we have
  \begin{equation}\label{eq:B0_Besov_pairing}
      \bigl|(f,g)_{\dot B^0_{2,2}(\sqsop)}\bigr| \le\|P_cf\|_{\dot B^{\rho_1-\sigma}_{r_1,2}(\sqsop)}\|P_cg\|_{\dot B^{\sigma-\rho_1}_{r_1',2}(\sqsop)}.
  \end{equation}
  Applying \eqref{eq:rho11}, or \eqref{eq:rho11_inf} when $q_1=\infty$, together with \eqref{eq:B0_Besov_pairing} and H\"older's inequality in time, we bound the right-hand side of \eqref{eq:energy_B0} by
  $$
    \lesssim \|P_cF\|_{L^{q_1'}(0,T;\dot B^{\sigma-\rho_1}_{r_1',2}(\sqsop))}^2 .
  $$
  Thus,
  \begin{equation}\label{eq:eq39}
    \left\|\int_0^T e^{-\mathrm{i}t\sqsop}
    P_cF(t)\,\tdt
    \right\|_{\dot{B}^{0}_{2,2}(\sqsop)}
    \lesssim
    \|P_cF\|_{L^{q_1'}(0,T;
    \dot{B}^{\sigma-\rho_1}_{r_1',2}(\sqsop))}.
  \end{equation}
  By duality and the preceding $L^2$ norm equivalence,
  $$
    \left\|e^{\mathrm{i}t\sqsop}P_cu_0
    \right\|_{L^{q_1}(0,T;
    \dot{B}^{\rho_1-\sigma}_{r_1,2}(\sqsop))}
    \lesssim
    \|P_cu_0\|_{\dot{B}^{0}_{2,2}(\sqsop)}.
  $$
  The case $q_1=\infty$ also follows directly from unitarity.
  The spectral lifting property
  $$
    \sop^\lambda:
    \dot{B}^{s}_{r,q}(\sqsop)
    \longrightarrow
    \dot{B}^{s-2\lambda}_{r,q}(\sqsop)
  $$
  now yields \eqref{eq:str4} by applying the preceding
  homogeneous estimate to $\sop^{\sigma/2}P_cu_0$.

  To obtain the energy-output estimate, apply
  \eqref{eq:eq39} on the fixed interval $(0,T)$ to the
  truncated source $\mathbf{1}_{(0,t)}F$.
  Unitarity then gives
  \begin{equation}\label{eq:rho13}
    \begin{aligned}
      \left\|\int_0^t e^{\mathrm{i}(t-\tau)\sqsop}P_cF(\tau)\,\td\tau\right\|_{L^\infty(0,T;\dot{B}^{0}_{2,2}(\sqsop))}
       &=\operatorname*{ess\,sup}_{0<t<T}\left\|\int_0^t e^{-\mathrm{i}\tau\sqsop}P_cF(\tau)\,\td\tau \right\|_{\dot{B}^{0}_{2,2}(\sqsop)}\\
      &\lesssim \|P_cF\|_{L^{q_1'}(0,T; \dot{B}^{\sigma-\rho_1}_{r_1',2}(\sqsop))}.
    \end{aligned}
  \end{equation}
  The constant is independent of the truncation time $t$.

  The homogeneous estimate and Minkowski's integral inequality
  also yield
  \begin{equation}\label{eq:rho12}
    \begin{aligned}
      &\left\|\int_0^t e^{\mathrm{i}(t-\tau)\sqsop}P_cF(\tau)\,\td\tau\right\|_{L^{q_1}(0,T;\dot{B}^{\rho_1-\sigma}_{r_1,2}(\sqsop))}\\
      \le&
      \int_0^T\left\|e^{\mathrm{i}(t-\tau)\sqsop}P_cF(\tau)\right\|_{L^{q_1}(\tau,T; \dot{B}^{\rho_1-\sigma}_{r_1,2}(\sqsop))}\,\td\tau \\
      \lesssim& \int_0^T \|P_cF(\tau)\|_{\dot{B}^{0}_{2,2}(\sqsop)} \,\td\tau   \\
      =&\|P_cF\|_{L^1(0,T;\dot{B}^{0}_{2,2}(\sqsop))}.
    \end{aligned}
  \end{equation}
  Here the homogeneous estimate is applied after translating $(\tau,T)$ to $(0,T-\tau)$ and restricting the estimate on $(0,T)$, so its constant is independent of $\tau$. For $q_1=\infty$, the same argument uses the essential supremum in time.

  We now interpolate these estimates in the spectral Besov scale. Recall that
  $$
    \rho_1-\sigma=-\frac2{q_1},
    \qquad
    1-\sigma-\rho_2=\frac2{q_2}.
  $$
  If $q_1\leq q_2$ and $q_1<\infty$, interpolation between \eqref{eq:rho11} and \eqref{eq:rho12}, with weight $q_1/q_2$ on the former, gives the desired source exponent while leaving the output space unchanged.
  If $q_2\leq q_1$ and $q_2<\infty$, we use the diagonal estimate and the energy-output estimate associated with the boundary pair $(q_2,r_2)$, with weight $q_2/q_1$ on the diagonal estimate. This leaves the source space unchanged and gives the desired output exponent. When $q_1=q_2=\infty$, the required bound is \eqref{eq:rho11_inf}. Hence, in all cases,
  $$
    \left\|\int_0^t e^{\mathrm{i}(t-\tau)\sqsop}
    P_cF(\tau)\,\td\tau
    \right\|_{L^{q_1}(0,T;
    \dot{B}^{\rho_1-\sigma}_{r_1,2}(\sqsop))}
    \lesssim
    \|P_cF\|_{L^{q_2'}(0,T;
    \dot{B}^{1-\sigma-\rho_2}_{r_2',2}(\sqsop))}.
  $$

  Applying this estimate to $\sop^{(\sigma-1)/2}P_cF$ and using spectral lifting gives
  $$
    \left\|\int_0^t
    \frac{e^{\mathrm{i}(t-\tau)\sqsop}}{\sqsop}
    P_cF(\tau)\,\td\tau
    \right\|_{L^{q_1}(0,T;
    \dot{B}^{\rho_1}_{r_1,2}(\sqsop))}
    \lesssim
    \|P_cF\|_{L^{q_2'}(0,T;
    \dot{B}^{-\rho_2}_{r_2',2}(\sqsop))}.
  $$
  All preceding arguments apply equally to $e^{-\mathrm{i}t\sqsop}P_c$, since the dispersive and energy estimates are unchanged. Taking the difference
  of the two half-wave estimates proves \eqref{eq:str5}.

  Finally, the cosine term in \eqref{eq:rep} follows from the two homogeneous half-wave estimates, and the initial velocity term follows by applying \eqref{eq:str4} to
  $\sop^{-1/2}P_cu_1$. Sobolev embedding gives the remaining admissible pairs, completing the proof. The calculations are first justified for spectrally localised data and time-simple source terms;
  the estimates extend to the stated spaces by density.
\end{proof}

By Lemma~\ref{lem_V}, if $\|V\|_{\K}<\pi$, then a fortiori $\|V\|_{\K}<4\pi$, so $\sop$ has purely absolutely continuous spectrum and no eigenvalues or resonances on $[0,+\infty)$,
and $P_c=\mathrm{Id}$. Since also $\|V_-\|_{\K}\le\|V\|_{\K}<\pi$, the hypotheses of Theorem \ref{eqv_besov} are satisfied, the estimate can be restated in the classical homogeneous Besov scale. In this case the estimate simplifies as follows.

\begin{cor}\label{cor_strichartz_smallness}
  Assume that
  $$
    V\in\K_0\cap L^{3/2}(\R^3), \qquad \|V\|_{\K}<\pi.
  $$
  Then, for every choice of indices $(q_i,r_i,\sigma,\rho_i)$ satisfying the conditions of Theorem~\ref{thm_strichartz},
  \begin{equation}\label{eq:str6}
    \begin{aligned}
      &\|u\|_{L^{q_1}(0,T;\dot{B}^{\rho_1}_{r_1,2}(\R^3))} +\|u\|_{L^{\infty}(0,T;\dot{H}^{\sigma}(\R^3))} +\|u_t\|_{L^{\infty}(0,T;\dot{H}^{\sigma-1}(\R^3))}   \\
      \lesssim&   \|u_0\|_{\dot{H}^{\sigma}(\R^3)}+\|u_1\|_{\dot{H}^{\sigma-1}(\R^3)}+\|F\|_{L^{q_2'}(0,T;\dot{B}^{-\rho_2}_{r_2',2}(\R^3))}.
    \end{aligned}
  \end{equation}
\end{cor}
The additional energy bounds in $L^\infty(0,T;\dot{H}^{\sigma}\times \dot{H}^{\sigma-1})$ follow by applying \eqref{eq:rho13} to the half-wave variables $W_\pm=P_cu_t\pm i\sqsop P_cu$, with $(q_1,r_1)$ replaced by
$(q_2,r_2)$ and with the same Sobolev reduction for non-boundary admissible pairs. Indeed, the source exponent before spectral lifting is $1-\sigma-\rho_2 =3(1/2-1/r_2)-1/q_2\in I_{r_2'}$
by \eqref{eq:wave_admissible}; after applying $\sop^{(\sigma-1)/2}$, the required source norm has index $-\rho_2$. The estimate holds for both signs, and recovering $P_cu_t$ and $\sqsop P_cu$ from $W_\pm$ gives the two $L^\infty_t$ terms.

By Remark~\ref{constant_smallV}, the Bernstein constants are uniform over $V\in\mathcal V$. Moreover, the small-Kato-norm dispersive
estimate for the Schr\"odinger propagator \cite[Theorem~2]{BG2012} gives
$$
  \|e^{it\sop}\varphi\|_{L^\infty(\R^3)}\le C\left(1-\frac{\kappa}{4\pi}\right)^{-2}|t|^{-3/2}\|\varphi\|_{L^1(\R^3)}.
$$
Thus the frequency-localized argument in \cite[Lemma~4.1]{bui} yields \eqref{eq:dispersive_estimate} with a constant uniform over $\mathcal V$. The proof of Theorem~\ref{thm_strichartz}, together
with the uniform Besov norm equivalence in Remark~\ref{constant_smallV}, then gives \eqref{eq:str6} with a constant depending only on $\kappa$, $M$, and the fixed indices and cutoffs.

\section{Well-posedness}\label{sec:wp}

This section is devoted to the well-posedness of weak solutions of \eqref{eq:problem}. We first recall the definition.
\begin{Def}\label{def_weak}
  A function $u(t)$ is a weak solution of \eqref{eq:problem} if $(u,u_t)\in L^\infty(0,T;\H)$ for any $T>0$, equation \eqref{eq:problem} is satisfied in the sense of distributions, i.e.,
  \begin{equation*}
    \int_{0}^{T}\int_{\R^3} \l(-u_t\phi_t+\ga u_t \phi + \nabla u \nabla \phi+mu\phi+f(u)\phi\r)\td x\td t=\int_{0}^{T}\la g,\phi\ra \tdt,
  \end{equation*}
  for any $\phi \in C_0^{\infty}((0,T)\times \R^3)$ and the initial data are attained strongly in the lower space, i.e.,
  $$
    \lim_{t\to0^+}\|(u(t)-u_0, u_t(t)-u_1)\|_{L^2(\R^3)\times H^{-1}(\R^3)}=0.
  $$
\end{Def}
\begin{remark}\label{rem:continuity}
  The regularity above provides a continuous representative in the lower space. The definition requires this representative to attain the prescribed initial data. Since $(u,u_t)\in L^\infty(0,T;\H)$, under the growth condition \eqref{F1}, we have $f(u) \in L^\infty(0,T;H^{-1}(\R^3))$. Thus, $u_{tt}\in L^\infty(0,T;H^{-1}(\R^3))$ and the equation holds in $L^\infty(0,T;H^{-1}(\R^3))$.
  Consequently, the time derivative of the phase variable satisfies
  $$
    (u_t, u_{tt}) \in L^\infty(0,T; L^2(\R^3) \times H^{-1}(\R^3)).
  $$
  Since $(u,u_t)\in W^{1,\infty}(0,T;L^2\times H^{-1}(\R^3))$, the vector-valued Sobolev embedding implies that $(u, u_{t})$ admits a continuous representative in the lower space:
  $$
    (u, u_t) \in C([0,T]; L^2(\R^3) \times H^{-1}(\R^3)).
  $$
  Combining this continuity with the energy bound $(u,u_t)\in L^\infty(0,T;\H)$, we obtain $(u,u_t)\in C_w([0,T];\H)$. The continuous representative in $L^2(\R^3)\times H^{-1}(\R^3)$ therefore has a trace at $t=0$.
  According to Definition~\ref{def_weak}, this trace equals the prescribed initial data:
  $$
    \lim_{t\to0^+}(u(t),u_t(t))=(u_0,u_1)
    \quad\text{strongly in }L^2(\R^3)\times H^{-1}(\R^3).
  $$
\end{remark}

The following lemma provides the stationary solution used in the decomposition.
\begin{lemma}\label{le_sta_sol}
  Let $f$ satisfy \eqref{F1} and \eqref{F2}. Then the elliptic equation
  \begin{equation}\label{eq:sta_sol}
    -\Delta h+ m h+f(h)=g\in H^{-1}(\R^3)
  \end{equation}
  has a solution $h\in H^1(\R^3)$, and any two solutions $h_1,h_2$ of \eqref{eq:sta_sol} satisfy $h_1-h_2\in L^q(\R^3)$ for every $2\le q<\infty$.
\end{lemma}
\begin{proof}
  By \eqref{F1}--\eqref{F2}, $0\le F(s)\le C(|s|^2+|s|^{p+1})$. Since $p+1<6$ and $m>0$, the energy functional is well defined
  and coercive on $H^1(\R^3)$. If $h_n\rightharpoonup h$ in $H^1(\R^3)$, choose a subsequence realizing the lower limit of $\int_{\R^3}F(h_n)\,\tdx$. Local Rellich compactness and a diagonal extraction give
  $h_{n_k}\to h$ almost everywhere. Since $F\ge0$, Fatou's lemma yields
  $$
    \int_{\R^3}F(h)\,\tdx
    \le\liminf_{n\to\infty}\int_{\R^3}F(h_n)\,\tdx.
  $$
  Together with weak lower semicontinuity of the quadratic terms and weak continuity of the forcing term, this proves weak lower semicontinuity of the energy functional. The direct method then gives a minimizer, whose Euler--Lagrange equation is \eqref{eq:sta_sol}.

  Let $h_1,h_2$ solve \eqref{eq:sta_sol} and set
  $$
    a(x):=-\int_0^1 f^{\prime}\bigl((1-\tau)h_1+\tau h_2\bigr)\,\td\tau ,
  $$
  so that $\bar{h}:=h_1-h_2$ satisfies $-\Delta \bar{h}+m\bar{h}=a\,\bar{h}$.

  Since $f\in C^3$, $f(0)=f^\prime(0)=0$ and $f(s)s\ge0$, Taylor expansion at the
  origin gives $f^{\prime\prime}(0)=0$. Hence,
  $$
    |f^\prime(s)|\le C\max_{[-1,1]}|f^{\prime\prime\prime}| |s|^2\quad(|s|\le1),
    \qquad
    |f^\prime(s)|\le C|s|^{p-1}\quad(|s|\ge1).
  $$
  Thus,
  \begin{equation}\label{eq:fprime_global}
    |f^\prime(s)|\le C\bigl(|s|^2+|s|^{p-1}\bigr),\qquad s\in\R .
  \end{equation}

  Since $3<p<5$, we have $3<12/(p-1)<6$, so $H^1(\R^3)\hookrightarrow L^{12/(p-1)}(\R^3)\cap L^6(\R^3)$. Applying \eqref{eq:fprime_global} to $s=(1-\tau)h_1+\tau h_2$ and using convexity,
  $$
    |a(x)|\le C\bigl(|h_1|^2+|h_2|^2+|h_1|^{p-1}+|h_2|^{p-1}\bigr),
  $$
  whence,
  $$
    \|a\|_{L^{6/(p-1)}(\R^3)}\le C\sum_{i=1}^2
    \bigl(\|h_i\|_{L^{12/(p-1)}(\R^3)}^2+\|h_i\|_{L^6(\R^3)}^{p-1}\bigr)<\infty .
  $$
  Since $6/(p-1)>3/2$, the standard de Giorgi--Moser iteration for $-\Delta \bar{h}+m\bar{h}=a\bar{h}$ gives $\bar{h}\in L^q(\R^3)$ for every $2\le q<\infty$.
\end{proof}
The second assertion shows that the space-time class of Theorem \ref{regular} is independent of the
stationary solution chosen: if $h_1$ and $h_2$ both solve \eqref{eq:sta_sol} then
$h_1-h_2\in L^{\frac{24}{5-p}}(\R^3)$, so $u-h_1$ and $u-h_2$ lie in the same space.

Let $h\in H^1(\R^3)$ be a stationary solution of
$$
  -\Delta h+mh+f(h)=g.
$$
Testing this equation with $h$ and using $f(s)s\ge0$ gives the following bound, uniform over all stationary solutions corresponding to the fixed force $g$:
\begin{equation}\label{eq:h_energy}
  \|h\|_{H^1(\R^3)} \le C\|g\|_{H^{-1}(\R^3)},
\end{equation}
where the constant $C$ only depends on $m$.

Set
\begin{equation}\label{eq:Vdef}
  V:=f_1^\prime(h).
\end{equation}
We have $\nabla V=f_1^{\prime\prime}(h)\nabla h$. From \eqref{eq:fdecomp_def} and \eqref{eq:fdecomp} one has $|f_1^{\prime\prime}(s)|\le C|s|^{p-2}$ for all $s\in\R$. Hence, by H\"older's inequality and $2(p-2)<6$,
\begin{equation}\label{eq:V_kato_small}
  \begin{aligned}
    \|\nabla V\|_{L^1(\R^3)}
     & \le C\int_{|h|\ge \Lambda/2}|h|^{p-2}|\nabla h|\td x      \le C\|\nabla h\|_{L^2(\R^3)}
    \left(\int_{|h|\ge \Lambda/2}|h|^{2(p-2)}\td x\right)^{1/2}                                \\
     & \le C\Lambda^{p-5}\|\nabla h\|_{L^2(\R^3)}
    \|h\|_{L^6(\R^3)}^3 \le C\|g\|_{H^{-1}(\R^3)}^4\,\Lambda^{p-5}.
  \end{aligned}
\end{equation}
Consequently,
$$\|V\|_{\K} \le \|\nabla V\|_{L^1(\R^3)}\le C\|g\|_{H^{-1}(\R^3)}^4\,\Lambda^{p-5}.$$
Choose $\Lambda=\Lambda(\|g\|_{H^{-1}(\R^3)})$ large enough so that
\begin{equation}\label{eq:N_choice}
  C\|g\|_{H^{-1}(\R^3)}^4\,\Lambda^{p-5}\le \frac{\pi}{2}.
\end{equation}
Then every stationary solution $h$ satisfies
\begin{equation}\label{eq:V-egy}
  \|V\|_{\K}=\|f_1^{\prime}(h)\|_{\K}\le \frac{\pi}{2}.
\end{equation}
Moreover,
$$
\|V\|_{L^{3/2}(\R^3)}=\|f_1'(h)\|_{L^{3/2}(\R^3)}\le C\|h\|_{L^{3(p-1)/2}}^{p-1}\le C\|h\|_{H^1}^{p-1}\le C\|g\|_{H^{-1}}^{p-1}.
$$
By Lemma~\ref{lem_V} and Remark~\ref{constant_smallV}, for fixed admissible indices $r,q,s$, the norm-equivalence constant in
Theorem~\ref{eqv_besov} can be chosen independently of $h$ and depends only on the fixed indices, the structural constants, and $\|g\|_{H^{-1}(\R^3)}$.

By Lemma~\ref{lem_V}, our choice $V=f_1^{\prime}(h)$ satisfies the hypotheses required for the Strichartz estimates \eqref{eq:str6}
to hold for $\sop=-\Delta+V$. We may therefore apply these estimates in the argument below.

For this choice, the Lipschitz constant of $f_2$ is bounded above by $C_{\Lambda}=C_2\Lambda^{p-1}$; see \eqref{eq:fdecomp}. Since
$\Lambda=\Lambda(\|g\|_{H^{-1}(\R^3)})$, all constants in the subsequent estimates may depend on $\Lambda$, hence on $\|g\|_{H^{-1}(\R^3)}$.

We shall use the energy functional
$$E(u)=\frac12\|u_t\|^2+\frac12\|\nabla u\|^2+\frac{m}{2}\|u\|^2+\int_{\R^3}F(u)\tdx-\la g,u\ra.$$
By \eqref{F2}, there exist constants $C, c>0$, such that
$$E(u)\ge C(\|u_t\|^2+\|\nabla u\|^2+\|u\|^2)-c(\|g\|^2_{H^{-1}(\R^3)}+1).$$

\subsection{The a priori estimate}

We now prove Theorem~\ref{regular}, that is, the space-time bound
\begin{equation*}
  \delta\in L^{\frac{8}{p-1}}(t_0,t_0+1;L^{\frac{24}{5-p}}(\R^3)),
\end{equation*}
where $\delta=u-h$ is the perturbation of a stationary solution $h$. This bound is used directly in the uniqueness argument and in the higher-regularity component of the asymptotic compactness proof.

\begin{proof}[Proof of Theorem \ref{regular}]
  Set
  $$R = \|(u(t),u_t(t))\|_{L^\infty(t_0,t_0+1;\H)} + \|g\|_{H^{-1}(\R^3)} + 1.$$

  The perturbation $\delta=u-h$ satisfies the equation
  \begin{equation}\label{eq:delta}
    \delta_{tt} -\Delta \delta+\ga \delta_t+m\delta+f(u) -f(h)=0, \,\delta(t_0)=u(t_0)-h, \, \delta_t (t_0) =u_t(t_0).
  \end{equation}
  The nonlinear term is treated by means of the decomposition $f = f_1 + f_2$ of \eqref{eq:fdecomp}, in which
  \begin{equation}
    | f^{\prime\prime\prime}_1( s)|\le C|s|^{p-3},\quad  f^{\prime\prime}_1(0)=0\quad\text{and}\quad | f^{\prime}_2(s)|\le C_{\Lambda}.
  \end{equation}
  By the fundamental theorem of calculus,
  \begin{equation}
    f(u) - f(h)  =f_1(u)-f_1(h)+f_2(u)-f_2(h)=  \delta \int_0^1 f_1^{\prime} ( h + \tau \delta) \td \tau+f_2(u)-f_2(h),
  \end{equation}
  and the integral term expands as
  \begin{equation}
    \begin{aligned}
      \int_0^1 f^{\prime}_1 (h + \tau \delta) \td \tau & =  \int_0^1 (f_1^{\prime} ( h + \tau \delta)-f^{\prime}_1(h)) \td \tau+f^{\prime}_1(h)                   \\
                                                       & =\delta \int_0^1 \int_0^\tau f^{\prime\prime}_1 ( h+\theta  \delta) \td \theta\td\tau  +f^{\prime}_1(h).
    \end{aligned}
  \end{equation}

  According to \eqref{eq:V-egy}, $\sop=-\Delta+V$ satisfies all the hypotheses of Corollary \ref{cor_strichartz_smallness}. Then \eqref{eq:delta} can be rewritten as
  \begin{equation}\label{eq:delta0}
    \delta_{tt} +\sop\delta= -\delta^2 \Gamma+f_2(h)-f_2(u)-\ga \delta_t-m\delta,
  \end{equation}
  where
  $$\Gamma=\int_0^1 \int_0^\tau f^{\prime\prime}_1 ( h+\theta  \delta) \td \theta\td\tau.$$
  Let $\delta_2$ be the unique solution of the linear wave equation
  \begin{equation}\label{eq:delta2}
    \ptt\delta_2+\sop\delta_2  = f_2(h)-f_2(u)-\ga \delta_t-m\delta,~~\delta_2(t_0)=u(t_0)-h, \, \partial_t\delta_2 (t_0) =u_t(t_0).
  \end{equation}
  Then $\delta_1=\delta-\delta_2$ satisfies a linear wave equation with vanishing datum
  \begin{equation}\label{eq:delta1}
    \ptt\delta_1+\sop\delta_1  = -\delta^2\Gamma, \quad \delta_1(t_0)=0, \, \partial_t\delta_1 (t_0) =0.
  \end{equation}
  Since $\delta_1\in L^\infty(t_0,t_0+1;H^1(\R^3))$, $\partial_t\delta_1\in L^\infty(t_0,t_0+1;L^2(\R^3))$, and
  $\delta^2\Gamma\in L^\infty(t_0,t_0+1;H^{-1}(\R^3))$, the standard linear representation applied to \eqref{eq:delta1} yields
  $$
    \delta_1(t) =-\int_{t_0}^t\frac{\sin((t-s)\sqsop)}{\sqsop}\,\delta(s)^2\Gamma(s)\,ds
  $$
  in $C([t_0,t_0+1];L^2(\R^3))$.
  We may therefore apply the inhomogeneous Strichartz estimate to this Duhamel solution whenever its forcing term belongs to the
  corresponding Besov space.

  The Strichartz estimate cannot be applied to achieve the target norm directly, because $f(u)$ need not have the required integrability. We therefore first establish the following intermediate bound
  \begin{equation*}
    \delta\in L^q(t_0,t_0+1;\dot B^0_{r,2}(\R^3)),
    \qquad
    \frac1q=\frac{p-3}{4},\quad \frac1r=\frac{5-p}{12},
  \end{equation*}
  and then reach the target by one further application of the estimate.

  The indices are chosen so that $\delta^2\Gamma$ is accessible to the product and composition rules in
  Besov spaces, subject to $q,r>2$, while $\tilde q$ is dictated by the admissibility condition
  \eqref{eq:wave_admissible}. The parameters $\rho_k$ form a finite ladder with a regularity gain of $\omega$ at each rung. This ladder starts at the value $\rho_0$ furnished by the Sobolev embedding and ends at the target $\rho_K=0$.

  Let $\omega=(5-p)/24$ and choose an integer $K$ with
  $$(K-1)\omega<\frac{p-3}{4}\le K\omega<\frac{p+1}{6},$$
  and set $\rho_k=(k-K)\omega$ for $k=0,1,\dots,K$, so that $-(p+1)/6<\rho_k\le 0$. Put
  $$
    \begin{aligned}
      \frac1r=\frac{5-p}{12},\quad \frac1{q}=\frac{p-3}{4},\quad \frac{1}{\tilde{{q}}}=\frac{p+1}{12},\quad \lambda_k=\rho_k+\frac{p+1}{6}\in (0,1).
    \end{aligned}
  $$

  The chain of Sobolev embeddings
  $$H^1(\R^3)=B^{1}_{2,2}(\R^3)\hookrightarrow B^{\rho_0+\frac{p+1}{4}}_{2,2}(\R^3)\hookrightarrow \dot{B}^{\rho_0+\frac{p+1}{4}}_{2,2}(\R^3)\hookrightarrow \dot{B}^{\rho_0}_{r,2}(\R^3)$$
  shows that $\delta \in L^{\infty}(0,T;\dot{B}^{\rho_0}_{r,2}(\R^3))$, which starts the ladder.

  Both $(q, r, \rho_k, \sigma_k)$ and $(\tilde{q}, r, -\lambda_k, 1-\sigma_k)$ satisfy the admissibility condition \eqref{eq:wave_admissible}, since
  $$\rho_k+3\l(\frac12-\frac1r\r)-\frac1{q}=\rho_k+1=\sigma_k=1-\l(-\lambda_k+3\l(\frac12-\frac1{r}\r)-\frac1{\tilde{q}}\r).$$

  Applying the Strichartz estimate \eqref{eq:str6} to \eqref{eq:delta1} gives
  \begin{equation}\label{eq:delta_str1}
    \|\delta_1\|_{L^{q}(t_0,t_0+1;\dot{B}^{\rho_k}_{r,2}(\R^3))}\lesssim \|\delta^2\Gamma\|_{L^{\tilde{q}'}(t_0,t_0+1;\dot{B}^{\lambda_k}_{r',2}(\R^3))}.
  \end{equation}
  The right-hand side is estimated by Proposition \ref{paraproduct}:
  \begin{equation}\label{eq:d2F1}
    \|\delta^2\Gamma\|_{\dot{B}^{\lambda_k}_{r',2}(\R^3)}\lesssim \|\delta^2\|_{\dot{B}^{\lambda_k}_{\frac{12}{11-p},2}(\R^3)}\|\Gamma\|_{\dot{B}^{\lambda_k}_{l_k,2}(\R^3)},
  \end{equation}
  where
  $$\frac{1}{l_k}=\frac{p-2}{6}+\frac{\lambda_k}{3}.$$
  With $s_1=(p+3)/8$, $s_2=\rho_{k-1}$ in Proposition \ref{negative_index}, we have $s_1+s_2=\lambda_k>0$, thereby controlling $\delta^2$ through a negative-order Besov norm,
  \begin{equation}
    \|\delta^2\|_{\dot{B}^{\lambda_k}_{\frac{12}{11-p},2}(\R^3)}\lesssim \|\delta\|_{\dot{B}^{\frac{3+p}{8}}_{2,\infty}(\R^3)}\|\delta\|_{\dot{B}^{\rho_{k-1}}_{r,2}(\R^3)}\lesssim R\|\delta\|_{\dot{B}^{\rho_{k-1}}_{r,2}(\R^3)}.
  \end{equation}
  Since $\Gamma=\int_0^1 \int_0^\tau f^{\prime\prime}_1 ( h+\theta  \delta) \td \theta\td\tau$, we have
  \begin{equation}
    \|\Gamma\|_{\dot{B}^{\lambda_k}_{l_k,2}(\R^3)}\le \int_0^1 \int_0^\tau \|f^{\prime\prime}_1 ( h+\theta  \delta) \|_{\dot{B}^{\lambda_k}_{l_k,2}(\R^3)}\td \theta\td \tau \le \sup_{0\le \theta\le 1}\|f^{\prime\prime}_1 ( h+\theta  \delta) \|_{\dot{B}^{\lambda_k}_{l_k,2}(\R^3)},
  \end{equation}
  while Proposition \ref{comp} gives
  \begin{equation}\label{eq:fpp}
    \begin{aligned}
      \|f^{\prime\prime}_1 ( h+\theta  \delta) \|_{\dot{B}^{\lambda_k}_{l_k,2}(\R^3)} & \lesssim  \|h+\theta  \delta\|_{\dot{B}^{\lambda_k}_{\frac{6}{1+2\lambda_k},2}(\R^3)}\||h+\theta  \delta|^{p-3}\|_{L^{\frac{6}{p-3}}(\R^3)} \\
                                                                                      & \lesssim \|h+\theta  \delta\|_{H^{1}(\R^3)}\||h+\theta  \delta|^{p-3}\|_{L^{\frac{6}{p-3}}(\R^3)}\lesssim R^{p-2}.
    \end{aligned}
  \end{equation}
  Combining \eqref{eq:d2F1}--\eqref{eq:fpp}, we arrive at
  \begin{equation}
    \|\delta^2\Gamma\|_{\dot{B}^{\lambda_k}_{r',2}(\R^3)}\lesssim R^{p-1}\|\delta\|_{\dot{B}^{\rho_{k-1}}_{r,2}(\R^3)}.
  \end{equation}
  Consequently,
  \begin{equation}
    \|\delta^2\Gamma\|_{L^{\tilde{q}'}(t_0,t_0+1;\dot{B}^{\lambda_k}_{r',2}(\R^3))}\lesssim R^{p-1}\|\delta\|_{L^{\tilde{q}'}(t_0,t_0+1;\dot{B}^{\rho_{k-1}}_{r,2}(\R^3))}.
  \end{equation}
  Since $q>\tilde{q}'$, this in turn gives
  \begin{equation}\label{eq:lambda_rho}
    \|\delta^2\Gamma\|_{L^{\tilde{q}'}(t_0,t_0+1;\dot{B}^{\lambda_k}_{r',2}(\R^3))}\lesssim R^{p-1}\|\delta\|_{L^{q}(t_0,t_0+1;\dot{B}^{\rho_{k-1}}_{r,2}(\R^3))}.
  \end{equation}
  Substituting \eqref{eq:lambda_rho} into \eqref{eq:delta_str1}, we obtain
  \begin{equation}
    \|\delta_1\|_{L^{q}(t_0,t_0+1;\dot{B}^{\rho_k}_{r,2}(\R^3))}\lesssim R^{p-1} \|\delta\|_{L^{q}(t_0,t_0+1;\dot{B}^{\rho_{k-1}}_{r,2}(\R^3))}.
  \end{equation}

  We turn to $\delta_2$. For $k=1,\dots,K$ let
  $$\frac{1}{q_k}=\frac1q+\rho_k .$$
  Then $q_k\ge q$, and $(q_k, r, \rho_k, 1)$ is admissible because
  $$\rho_k+3\l(\frac12-\frac1r\r)-\frac{1}{q_k}=3\l(\frac12-\frac1r\r)-\frac{1}{q}=1.$$
  The Strichartz estimate \eqref{eq:str6} applied to \eqref{eq:delta2} yields
  \begin{equation}\label{eq:delta_str2}
    \|\delta_2\|_{L^{q_k}(t_0,t_0+1;\dot{B}^{\rho_k}_{r,2}(\R^3))}\lesssim \|(u(t_0)-h, u_t(t_0))\|_{\H}+\|f_2(h)-f_2(u)-\ga \delta_t-m\delta\|_{L^1(t_0,t_0+1;L^2(\R^3))}.
  \end{equation}
  For the source term on the right of \eqref{eq:delta_str2},
  \begin{equation}
    \|\ga \delta_t+m\delta\|_{L^1(t_0,t_0+1;L^2(\R^3))}\lesssim  R,
  \end{equation}
  and
  \begin{equation}\label{eq:eqF2}
    \|f_2(h)-f_2(u)\|_{L^1(t_0,t_0+1;L^2(\R^3))} \le C_\Lambda R.
  \end{equation}
  Hence,
  \begin{equation}
    \|\delta_2\|_{L^{q_k}(t_0,t_0+1;\dot{B}^{\rho_k}_{r,2}(\R^3))}\le C_\Lambda R,
  \end{equation}
  and, since $q_k\ge q$ on a time interval of length one,
  \begin{equation}
    \|\delta_2\|_{L^{q}(t_0,t_0+1;\dot{B}^{\rho_k}_{r,2}(\R^3))}\le \|\delta_2\|_{L^{q_k}(t_0,t_0+1;\dot{B}^{\rho_k}_{r,2}(\R^3))}\le C_\Lambda R.
  \end{equation}

  The two estimates together give the recursion
  \begin{equation}\label{eq:iter}
    \|\delta\|_{L^{q}(t_0,t_0+1;\dot{B}^{\rho_k}_{r,2}(\R^3))}\lesssim R^{p-1} \|\delta\|_{L^{q}(t_0,t_0+1;\dot{B}^{\rho_{k-1}}_{r,2}(\R^3))}+C_\Lambda R ,
  \end{equation}
  which, iterated $K$ times from the initial bound at $\rho_0$, yields
  \begin{equation}
    \|\delta\|_{L^{q}(t_0,t_0+1;\dot{B}^{\rho_K}_{r,2}(\R^3))}\lesssim R^{(p-1)K+1}+C_\Lambda R^{(p-1)K+1}\lesssim C_\Lambda R^{(p-1)K+1}.
  \end{equation}
  As $r = \frac{12}{5-p} > 2$, the embedding $\dot{B}^{\rho_K}_{r,2}(\R^3)=\dot{B}^{0}_{r,2}(\R^3) \hookrightarrow L^r(\R^3)$ gives
  \begin{equation}\label{eq:str0}
    \|\delta\|_{L^{q}(t_0,t_0+1;L^r(\R^3))}\lesssim C_\Lambda R^{(p-1)K+1}.
  \end{equation}

  We emphasise that no a priori Besov regularity is assumed for $\delta$. Starting from the energy bound, $\delta\in L^\infty(t_0,t_0+1;\dot B^{\rho_0}_{r,2}(\R^3))$, and each subsequent Besov bound follows inductively from the product and composition estimates and the linear Strichartz estimates applied to \eqref{eq:delta2} and \eqref{eq:delta1}. Thus, the finite iteration derives all the Besov regularity used above directly from the weak formulation and the energy bound.

  In order to achieve the estimate in the norm $L^{\frac8{p-1}}(t_0,t_0+1;L^{\frac{24}{5-p}}(\R^3))$, we apply the Strichartz estimate \eqref{eq:str6} to equation \eqref{eq:delta0} with the admissible pair
  $$\l(\frac{8}{p-1},\quad \frac{24}{5-p}\r).$$
  This yields
  \begin{equation}
    \begin{aligned}
               & \|\delta\|_{L^{\frac{8}{p-1}}(t_0,t_0+1; L^{\frac{24}{5-p}}(\R^3))}                                                  \\
      \lesssim & \| (u(t_0)-h, u_t(t_0))\|_{\H} + \|\delta^2 \Gamma+f_2(u)-f_2(h)+\ga \delta_t+m\delta \|_{L^1(t_0,t_0+1;L^2(\R^3))}.
    \end{aligned}
  \end{equation}

  By \eqref{eq:fdecomp} we have $|\Gamma| \lesssim |h|^{p-2} + |u|^{p-2}$, so H\"older's inequality gives
  $$
    \|\delta^2 \Gamma\|_{L^2(\R^3)} \le \|\Gamma\|_{L^{\frac{6}{p-2}}(\R^3)}\|\delta\|^2_{L^{r}(\R^3)}\le CR^{p-2}\|\delta\|^2_{L^{r}(\R^3)},
  $$
  and integration in time yields
  \begin{equation}\label{eq:gammaL2}
    \|\delta^2\Gamma\|_{L^1(t_0,t_0+1;L^2(\R^3))}  \le C R^{p-2}  \|\delta\|_{L^{q}(t_0,t_0+1;L^{r}(\R^3))}^{2}.
  \end{equation}
  The remaining terms satisfy
  \begin{equation}\label{eq:f2L2}
    \|f_2(u)-f_2(h)+\ga \delta_t+m\delta\|_{L^1(t_0,t_0+1;L^2(\R^3))}\le C_\Lambda R.
  \end{equation}
  Hence,
  \begin{equation}
    \|\delta\|_{L^{\frac{8}{p-1}}(t_0,t_0+1; L^{\frac{24}{5-p}}(\R^3))}\lesssim  R^{p-2}\|\delta\|^2_{L^{q}(t_0,t_0+1;L^r(\R^3))}+C_\Lambda R\lesssim C_\Lambda^2 R^{2(p-1)K+p}.
  \end{equation}

  Finally, we show that the a priori estimate implies the strong time continuity in the energy space $\H$.
  The bounds in \eqref{eq:regu_est} and \eqref{eq:gammaL2} ensure that the terms on the right-hand side of equation \eqref{eq:delta0} satisfy:
  $$
    \delta^2 \Gamma \in L^1(0,T;L^2(\R^3)), \quad f_2(u)-f_2(h) \in L^\infty(0,T;L^2(\R^3)), \quad \delta_t \in L^\infty(0,T;L^2(\R^3)).
  $$
  Thus, the source term $\G(s) := -\gamma \delta_t -m\delta - \delta^2\Gamma - (f_2(u)-f_2(h))$ belongs to $L^1(0,T; L^2(\R^3))$.
  The perturbation $\delta$ solves the linear equation $\delta_{tt} + \sop \delta = \G$ with initial data $(u_0-h, u_1) \in \H$.
  By the Duhamel principle, we have:
  $$
    \delta(t) = \cos(t\sqsop)(u_0-h) + \frac{\sin(t\sqsop)}{\sqsop}u_1 + \int_{0}^t \frac{\sin((t-s)\sqsop)}{\sqsop} \G(s) \td s,
  $$
  and
  $$
    \delta_t(t) = -\sqsop\sin(t\sqsop)(u_0-h) + \cos(t\sqsop)u_1 + \int_{0}^t \cos((t-s)\sqsop) \G(s) \td s.
  $$
  We work with the homogeneous energy-space realization in which $\sqsop$ induces an isomorphism $\sqsop:\dot H^1(\R^3)\longrightarrow L^2(\R^3)$. By the spectral theorem, the operators
  $$
    \cos(t\sqsop),\qquad \frac{\sin(t\sqsop)}{\sqsop}
  $$
  are well-defined on the corresponding homogeneous Sobolev scale. In
  particular, the wave propagator
  $$
    \mathcal U(t)=\begin{pmatrix}
      \cos(t\sqsop)        & \quad\dfrac{\sin(t\sqsop)}{\sqsop} \\[2ex]
      -\sqsop\sin(t\sqsop) & \quad\cos(t\sqsop)
    \end{pmatrix}
  $$
  defines a strongly continuous group on $\dot H^1(\R^3)\times L^2(\R^3)$. Consequently, the homogeneous terms belong to $C([0,T];\dot H^1(\R^3)\times L^2(\R^3))$.

  For the inhomogeneous terms, because
  $$
    \|\sin(\tau \sqsop)\|_{L^2\to L^2}\le 1,
    \qquad
    \|\cos(\tau \sqsop)\|_{L^2\to L^2}\le 1,
  $$
  and both operator families are strongly continuous on $L^2(\R^3)$, dominated convergence implies that
  $$t \mapsto \int_{0}^t \frac{\sin((t-s)\sqsop)}{\sqsop} \G(s) \td s\quad \text{and}\quad t \mapsto \int_{0}^t \cos((t-s)\sqsop) \G(s) \td s$$
  define continuous maps into $\dot H^1(\R^3)$ and $L^2(\R^3)$, respectively.
  Consequently,
  $$(\delta, \delta_t) \in C([0,T]; \dot{H}^1(\R^3)\times L^2(\R^3)).$$
  Combining this with $\delta(t)=\delta(0)+\int_0^t\delta_t(\tau)\td \tau$ gives $(\delta, \delta_t) \in C([0,T]; \H)$. Since $h$ is time-independent, we conclude that $(u,u_t) \in C([0,T]; \H)$.
  This completes the proof.
\end{proof}
\begin{remark}\label{rem_iteration}
  Two obstructions appear as $p\to5$. The first is that the ladder has nowhere to stand: its rungs
  are measured in $L^{q}(t_0,t_0+1;\dot B^{\rho_k}_{r,2}(\R^3))$ with
  $$
    \frac1q=\frac{p-3}{4},\qquad \frac1r=\frac{5-p}{12},
  $$
  so that $(q,r)\to(2,\infty)$, the excluded endpoint for the three-dimensional wave Strichartz estimate.

  The second is that the ladder has no rise. At $p=5$ the increment $\omega$ vanishes, so the source and target indices coincide and the same norm stands on both sides of \eqref{eq:iter}. The recursion
  becomes self-referential and can no longer raise regularity. A contraction argument, when available, is instead applied to the difference $\delta=u-v$ of two solutions with the same initial data. Contraction inequalities of this type have been used to establish uniqueness at the critical exponent in higher dimensions; see Planchon \cite{planchon03} for
  $d=4,5$ and Bulut et al. \cite{bulut} for $d\ge6$.
\end{remark}

\subsection{Uniqueness}
The a priori space-time estimate \eqref{eq:regu_est} now yields continuous dependence on the initial data and uniqueness of weak solutions.
\begin{thm}\label{thm_uniq}
  Under the hypotheses of Theorem \ref{regular}, let $u$ and $v$ be weak solutions with initial data $(u_0,u_1)$ and $(v_0,v_1)$ in $\H$, and put $\bar u=u-v$. Then
  \begin{equation}\label{eq:con_depd}
    \|(\bar{u}(t),\bar{u}_t(t))\|_{\H}\le C_T \| (\bar{u}(0),\bar{u}_t(0))\|_{\H},\quad t\in [0,T],
  \end{equation}
  where $C_T$ depends on $T$, the $C([0,T];\H)$ norms of $(u,u_t)$ and $(v,v_t)$, and on $g$ through the cut-off parameter $\Lambda=\Lambda(\|g\|_{H^{-1}(\R^3)})$.
\end{thm}
Uniqueness here is unconditional: the two solutions are arbitrary weak solutions in the sense of
Definition \ref{def_weak}, and the space-time regularity used in the proof is not assumed but
supplied by Theorem \ref{regular}. See Remark \ref{rem_iteration} for the endpoint obstruction
that the present iteration encounters at $p=5$.

\begin{proof}
  Set
  $$
    \begin{aligned}
      R = & \sup_{t \in [0,T]} \left\{ \|(u(t),u_t(t))\|_{\H}, \|(v(t),v_t(t))\|_{\H} \right\}+\|h\|_{H^1(\R^3)}            \\
          & + \|u - h\|_{L^{\frac{8}{p-1}}(0,T;L^{\frac{24}{5-p}})}+ \|v - h\|_{L^{\frac{8}{p-1}}(0,T;L^{\frac{24}{5-p}})}.
    \end{aligned}
  $$

  The difference $\bar u=u-v$ satisfies
  \begin{equation}\label{eq:difference}
    \begin{cases}
      \bar{u}_{tt} + \ga \bar{u}_t - \Delta \bar{u} + m \bar{u} + f(u) - f(v) = 0, \\
      \bar{u}(0) = u_0 - v_0, \quad \bar{u}_t(0) = u_1 - v_1.
    \end{cases}
  \end{equation}

  We again split $f = f_1 + f_2$ and expand $f_1$ around the stationary solution $h$. The $f_2$ term is handled using its global Lipschitz bound. By the fundamental theorem of calculus,
  \begin{equation}
    \begin{aligned}
      f_1(u)-f_1(v)
       & = \bar{u} \int_0^1 f_1^{\prime}(v + \tau \bar{u})\td\tau                                                           \\
       & = f_1^{\prime}(h) \bar{u}+\bar{u} \int_0^1 \left[ f_1^{\prime}(v + \tau \bar{u}) - f_1^{\prime}(h) \right] \td\tau \\
       & = V \bar{u} + \bar{u}R_1(u,v;h),
    \end{aligned}
  \end{equation}
  where $V = f_1^{\prime}(h)$ and
  \begin{equation}\label{eq:R1_def}
    R_1(u,v;h):=\int_0^1\int_0^1f_1^{\prime\prime}\l(h+\theta(v - h+\tau \bar u)\r)(v - h+\tau \bar u)\,
    \td\theta\,\td\tau .
  \end{equation}
  For $f_2$ we write
  \begin{equation}
    f_2(u)-f_2(v)= \bar{u} \int_0^1 f_2^{\prime}(v + \tau \bar{u})\td\tau=:\bar u F_2,
  \end{equation}
  and $\|F_2\|_{L^\infty}\le C_{\Lambda}$ by \eqref{eq:fdecomp}.

  Then equation \eqref{eq:difference} becomes
  \begin{equation}\label{eq:decomposed}
    \begin{cases}
      \ptt \bar{u} - \Delta \bar{u}  + V \bar{u} = -\bar{u} R_1(u,v;h) - (m+F_2) \bar{u}- \ga \bar{u}_t, \\
      \bar{u}(0) = u_0 - v_0, \quad \bar{u}_t(0) = u_1 - v_1.
    \end{cases}
  \end{equation}
  The Strichartz estimate \eqref{eq:str6} applied to \eqref{eq:decomposed} on $[0,t]$ yields
  \begin{equation}\label{eq:strichartz}
    \begin{aligned}
               & \|\bar{u}\|_{L^{\frac{8}{p-1}}(0,t; L^{\frac{24}{5-p}}(\R^3))}+\|\nabla \bar{u}(t)\| + \|\bar{u}_t(t)\|               \\
      \lesssim & \| (\bar{u}(0),\bar{u}_t(0))\|_{\H} + \| \bar{u} R_1(u,v;h) + (m+F_2) \bar{u} + \ga \bar{u}_t\|_{L^1(0,t;L^2(\R^3))}.
    \end{aligned}
  \end{equation}
  By the growth condition on $f$, we have
  $$|R_1(u,v;h)| \lesssim \l(|h| + |v| + |u|\r)^{p-2}\l(|v-h|+|u-h|\r).$$
  By H\"older's inequality,
  $$\begin{aligned}
      \|R_1(u,v;h)\|_{L^{\frac{24}{p+7}}(\R^3)} & \lesssim \left\|\l(|h| + |v| + |u|\r)^{p-2}\right\|_{L^{\frac{12}{p+1}}(\R^3)} \l(\|v-h\|_{L^{\frac{24}{5-p}}(\R^3)}+\|u-h\|_{L^{\frac{24}{5-p}}(\R^3)}\r) \\
                                                & \lesssim R^{p-2}\l(\|v-h\|_{L^{\frac{24}{5-p}}(\R^3)}+\|u-h\|_{L^{\frac{24}{5-p}}(\R^3)}\r).
    \end{aligned}$$
  Thus,
  \begin{equation}
    \begin{aligned}
      \|\bar{u} R_1(u,v;h)\| & \lesssim  \|R_1(u,v;h)\|_{L^{\frac{24}{p+7}}(\R^3)}\|\bar{u}\|_{L^{\frac{24}{5-p}}(\R^3)}                                           \\
                             & \lesssim  R^{p-2}\l(\|v-h\|_{L^{\frac{24}{5-p}}(\R^3)}+\|u-h\|_{L^{\frac{24}{5-p}}(\R^3)}\r)\|\bar{u}\|_{L^{\frac{24}{5-p}}(\R^3)}.
    \end{aligned}
  \end{equation}
  Since $q=8/(p-1)$ and $1-2/q=(5-p)/4$, another H\"older inequality in time yields
  \begin{equation}\label{eq:nonlinear1}
    \begin{aligned}
               & \| \bar{u} R_1(u,v;h) \|_{L^1(0,t;L^2(\R^3))}                                                                                                                                                                        \\
      \lesssim & R^{p-2} t^{\frac{5-p}{4}} \l(\|v - h\|_{L^{\frac{8}{p-1}}(0,t;L^{\frac{24}{5-p}}(\R^3))}+\|u - h\|_{L^{\frac{8}{p-1}}(0,t;L^{\frac{24}{5-p}}(\R^3))}\r)\|\bar{u}\|_{L^{\frac{8}{p-1}}(0,t;L^{\frac{24}{5-p}}(\R^3))} \\
      \lesssim & R^{p-1} t^{\frac{5-p}{4}} \|\bar{u}\|_{L^{\frac{8}{p-1}}(0,t;L^{\frac{24}{5-p}}(\R^3))}.
    \end{aligned}
  \end{equation}
  As for the potential term $(m+F_2) \bar{u}$,
  \begin{equation}\label{eq:potential2}
    \|(m+F_2) \bar{u}\|_{L^1(0,t;L^2(\R^3))} \lesssim \|m+F_2\|_{L^\infty} \|\bar{u}\|_{L^1(0,t;L^2(\R^3))} \le C C_{\Lambda} \|\bar{u}\|_{L^1(0,t;H^1(\R^3))}.
  \end{equation}
  Now choose $S\in (0, 1)$ sufficiently small, depending only on $R$, $p$ and the constants in \eqref{eq:nonlinear1}, so that
  $$C R^{p-1} S^{\frac{5-p}{4}}\le \frac12.$$
  Then, for every $t\in[0,S]$, the factor $t^{\frac{5-p}{4}}\le S^{\frac{5-p}{4}}$ in \eqref{eq:nonlinear1} ensures that the term
  $\bar u\,R_1(u,v;h)$ can be absorbed into the left-hand side of \eqref{eq:strichartz}.
  Consequently,
  \begin{equation}\label{eq:combined2}
    \|\nabla \bar{u}(t)\| + \|\bar{u}_t(t)\|  \lesssim  \| (\bar{u}(0),\bar{u}_t(0))\|_{\H}+\| (m+F_2) \bar{u} + \ga \bar{u}_t \|_{L^1(0,t;L^2(\R^3))}.
  \end{equation}
  Substituting \eqref{eq:potential2} into \eqref{eq:combined2}, we obtain
  \begin{equation*}
    \|\nabla \bar{u}(t)\| + \|\bar{u}_t(t)\| \le C \l( \| (\bar{u}(0),\bar{u}_t(0))\|_{\H} + C_{\Lambda}\|(\bar{u},\bar{u}_t)\|_{L^1(0,t;\H)}\r),\quad 0 \le t \le S.
  \end{equation*}
  Combining this with the $L^2$ estimate
  $$\|\bar{u}(t)\|\le \|\bar{u}(0)\|+\int_0^t\|\bar{u}_t(s)\|\td s,$$
  we deduce that, for any $t\in [0, S]$,
  \begin{equation}\label{eq:final-energy}
    \|(\bar{u}(t),\bar{u}_t(t))\|_{\H} \le  C \l( \| (\bar{u}(0),\bar{u}_t(0))\|_{\H} + C_{\Lambda}\|(\bar{u},\bar{u}_t)\|_{L^1(0,t;\H)}\r).
  \end{equation}
  Applying Gronwall's inequality to \eqref{eq:final-energy} gives
  \begin{equation}\label{eq:gronwall}
    \|(\bar{u}(t),\bar{u}_t(t))\|_{\H}  \le C e^{C C_{\Lambda} t} \| (\bar{u}(0),\bar{u}_t(0))\|_{\H}, \quad 0 \le t \le S.
  \end{equation}
  In particular, there exists $C_0=C_0(\|g\|_{H^{-1}(\R^3)}) > 0$ such that
  $$
    \|(\bar{u}(t),\bar{u}_t(t))\|_{\H} \le C_0 \| (\bar{u}(0),\bar{u}_t(0))\|_{\H}, \quad 0 \le t \le S.
  $$

  The estimate extends to $[0,T]$ by iteration. Since $S$ depends only on $R$, $p$ and the constant in \eqref{eq:nonlinear1}, the same argument applied to the solutions shifted in time gives, for every $s \in [0,T-S]$,
  $$
    \|(\bar{u}(t),\bar{u}_t(t))\|_{\H}  \le C_0 \| (\bar{u}(s),\bar{u}_t(s))\|_{\H}, \quad s \le t \le s+S,
  $$
  and after $k = \lceil T/S \rceil$ steps
  $$
    \|(\bar{u}(t),\bar{u}_t(t))\|_{\H}  \le C_0^k \| (\bar{u}(0),\bar{u}_t(0))\|_{\H}, \quad 0 \le t \le T.
  $$
  The shifted solutions satisfy the same a priori bound on every interval of length $S$; hence the constants in the local estimate are independent of the shift time. Taking $C_T=C_0^{\lceil T/S\rceil}$ proves \eqref{eq:con_depd}.
\end{proof}

\subsection{Existence of weak solutions}
We construct weak solutions by approximation, in the spirit of \cite[Theorem 3.4]{lcc1}.
\begin{thm}[Existence]\label{exist}
  Let $g\in H^{-1}(\R^3)$ and let the nonlinearity $f$ satisfy \eqref{F1} and \eqref{F2}. Then for every $(u_0,u_1) \in \H$ there is a global weak solution $u$ of \eqref{eq:problem}, and
  \begin{equation}\label{eq:eng_est}
    \| (u,u_t)\|_{\H}\le Q(\| (u_0,u_1)\|_{\H}+\|g\|_{H^{-1}(\R^3)}),
    \qquad t\ge 0,
  \end{equation}
  where $Q$ is a continuous increasing function.
\end{thm}
\begin{proof}
  Choose regularised data $\{(u^n_0,u^n_1)\}\subset H^2(\R^3)\times H^1(\R^3)$ and $\{g^n\}\subset L^2(\R^3)$ with $(u^n_0,u^n_1)\to (u_0,u_1)$ in $\H$ and $g^n\to g$ in $H^{-1}(\R^3)$. By \cite[Proposition 1.1]{ef} the approximating problem
  \begin{equation}\label{eq:approx}
    \begin{cases}
      u^n_{tt}+\ga u^n_t-\Delta u^n+m u^n+f(u^n)=g^n, \\
      u^n(0)=u^n_0,u^n_t(0)=u^n_1,
    \end{cases}
  \end{equation}
  has a unique strong solution $u^n\in L^\infty_{loc}(\R^{+};H^2(\R^3))$.

  Define the approximate energy by
  $$
    \mathcal E_n(t):= \frac12\|u_t^n(t)\|^2+\frac12\|\nabla u^n(t)\|^2+\frac m2\|u^n(t)\|^2+\int_{\R^3}F(u^n(t,x))\,dx-\la g^n,u^n(t)\ra,
  $$
  where $F(s)=\int_0^s f(\tau)\,d\tau$. Since $g^n\in L^2(\R^3)\subset H^{-1}(\R^3)$, the forcing term is included in the potential energy. The energy identity is
  $$
    \frac{d}{dt}\mathcal E_n(t)=-\ga\|u_t^n(t)\|^2.
  $$
  By \eqref{F2}, $F(s)\ge0$. Moreover, the sequence $\{g^n\}$ is bounded in $H^{-1}(\R^3)$, and
  $$
    |\la g^n,u^n\ra|\le \varepsilon\|u^n\|_{H^1(\R^3)}^2+C_\varepsilon\|g^n\|_{H^{-1}(\R^3)}^2.
  $$
  Consequently, the energy identity and the uniform bounds on the initial data and on $g^n$ yield
  \begin{equation}\label{eq:approx_energy}
    \|(u^n(t),u_t^n(t))\|_{\H} \le Q_0\bigl(\|(u_0,u_1)\|_{\H}+\|g\|_{H^{-1}(\R^3)}\bigr),  \qquad t\ge0,
  \end{equation}
  where $Q_0$ is positive and increasing, and is independent of $n$ and $t$. In particular, no uniform bound on $\|g^n\|$ is required.

  By \eqref{eq:approx_energy} and the approximate equation \eqref{eq:approx}, for every $T>0$ there is a constant $C_T$, independent of $n$, such that
  $$
    \sup_n\|u_t^n\|_{L^\infty(0,T;L^2(\R^3))} +\sup_n\|u_{tt}^n\|_{L^\infty(0,T;H^{-1}(\R^3))}  \le C_T.
  $$
  Indeed, $\Delta u^n$, $u^n$, $f(u^n)$, and $g^n$ are uniformly bounded in $L^\infty(0,T;H^{-1}(\R^3))$; for $f(u^n)$ this follows from the growth
  assumptions and the uniform $H^1$ bound. Thus, for all $s,t\in[0,T]$,
  $$
    \|u^n(t)-u^n(s)\|+\|u_t^n(t)-u_t^n(s)\|_{H^{-1}(\R^3)}\le C_T|t-s|.
  $$
  In particular,
  $$
    \|u^n(t)-u_0^n\|+\|u_t^n(t)-u_1^n\|_{H^{-1}(\R^3)} \le C_Tt,\qquad 0\le t\le T.
  $$

  For each $k\in\mathbb N$, the compact embeddings
  $$
    H^1(B_k)\hookrightarrow\hookrightarrow L^2(B_k) \hookrightarrow H^{-1}(B_k)
  $$
  and the uniform bound \eqref{eq:approx_energy} allow us, by the Aubin--Lions--Simon theorem, to pass to a subsequence
  such that $u^n\to u$ in $C([0,T];L^2(B_k))$ for every $k$. A diagonal extraction and interpolation with the uniform
  $L^\infty(0,T;L^6(\R^3))$ bound then give
  $$
    u^n\to u\quad\text{in }L^q_{\mathrm{loc}}((0,T)\times\R^3) \quad\text{for every }q<6,
  $$
  and almost everywhere on $(0,T)\times\R^3$. In addition, after passing to a further subsequence,
  \begin{align*}
     & u^n\to u \text{ weakly star in } L^\infty(0,T;H^1(\R^3)),     \\
     & u^n_t\to u_t \text{ weakly star in } L^\infty(0,T;L^2(\R^3)).
  \end{align*}
  Since $p<5$, the strong local convergence in $L^5$ and the growth assumption on $f$ imply
  $$
    f(u^n)\to f(u) \quad\text{in }L^1_{\mathrm{loc}}((0,T)\times\R^3).
  $$
  Passing to the limit in the weak formulation of \eqref{eq:approx} shows that $u$ solves \eqref{eq:problem}. Weak lower semicontinuity
  gives the essential-supremum energy bound for the limit.

  It remains to verify the initial trace in the global lower space. The uniform time-Lipschitz estimates above, together with
  $(u_0^n,u_1^n)\to(u_0,u_1)$ in $H^1(\R^3)\times L^2(\R^3)$, imply
  $$
    \|u^n(t)-u_0\|+\|u_t^n(t)-u_1\|_{H^{-1}(\R^3)} \le C_Tt +\|u_0^n-u_0\| +\|u_1^n-u_1\|_{H^{-1}(\R^3)}.
  $$

  The uniform bounds and equi-Lipschitz estimates yield, by diagonal extraction on a countable dense subset of $[0,T]$ and extension by
  equi-Lipschitz continuity, a subsequence converging weakly in $L^2(\R^3)\times H^{-1}(\R^3)$ at every $t\in[0,T]$, with a weakly
  continuous limit. The convergence $u^n\to u$ in $C([0,T];L^2(B_k))$ identifies its first component with $u(t)$ for every $t$. 
  The weak-star convergence of $u_t^n$ identifies its second component with $u_t$ almost everywhere in time; continuity into $H^{-1}(\R^3)$
  extends this equality to every $t$.

  Weak lower semicontinuity now gives
  $$
    \|u(t)-u_0\| +\|u_t(t)-u_1\|_{H^{-1}(\R^3)} \le C_Tt,\qquad 0\le t\le T.
  $$
  Hence, $(u(t),u_t(t))\to(u_0,u_1)$ strongly in $L^2(\R^3)\times H^{-1}(\R^3)$ as $t\to0^+$, so the limit satisfies Definition~\ref{def_weak}.
  Applying Theorem~\ref{regular} gives $(u,u_t)\in C([0,T];\H)$; thus the essential-supremum estimate extends to every $t\in[0,T]$, proving \eqref{eq:eng_est}.
\end{proof}

The well-posedness established above allows us to define the solution semigroup $S(t):\H\rightarrow \H$ of \eqref{eq:problem} by
\begin{eqnarray*}
  S(t)(u_0,u_1)=(u(t),u_t(t)),
\end{eqnarray*}
where $u$ is the unique weak solution with data $(u_0,u_1)$. By Theorems~\ref{thm_uniq} and~\ref{exist}, $S(t)$ is well-defined and continuous on $\H$.

\section{Existence of the global attractor}\label{sec:attractor}
In this section we prove the existence of the global attractor for \eqref{eq:problem}, assuming \eqref{F1}--\eqref{F3} on the nonlinearity and $g\in H^{-1}(\R^3)$. The proof rests on the dissipativity of the semigroup, uniform tail estimates, and a refined decomposition of the solutions supplying the extra compactness.

\subsection{The weighted energy equality}
We establish the weighted energy identity by time-mollification arguments.
\begin{lemma}\label{le_energy_eqv}
  Let $g\in H^{-1}(\R^3)$, let the nonlinearity $f$ satisfy \eqref{F1} and \eqref{F2}, and let $u$ be a weak solution of
  \eqref{eq:problem} on $[0,T]$. Then for every nonnegative weight $\rho\in C^{\infty}(\R^3)$ with $\rho,\nabla\rho\in L^\infty(\R^3)$, the weighted energy functional
  $$
    E_\rho(u(t))
    :=\frac12\int_{\R^3}\rho(x)
    \left(|u_t|^2+|\nabla u|^2+m|u|^2+2F(u)\right)\td x
    -\la g,\rho u(t)\ra
  $$
  is absolutely continuous on $[0,T]$, and, for every
  $0\le s\le t\le T$,
  \begin{equation}\label{eq:energy_eqv_int}
    E_\rho(u(t))-E_\rho(u(s))
    +\ga\int_s^t\int_{\R^3}\rho|u_\tau|^2\td x\td\tau
    +\int_s^t\int_{\R^3}
    (\nabla u\cdot\nabla\rho)u_\tau\td x\td\tau
    =0 .
  \end{equation}
  Consequently,
  \begin{equation}\label{eq:energy_eqv}
    \frac{\td}{\td t}E_\rho(u)
    +\ga\int_{\R^3}\rho|u_t|^2\td x
    +\int_{\R^3}(\nabla u\cdot\nabla\rho)u_t\td x
    =0
  \end{equation}
  for almost every $t\in[0,T]$.
\end{lemma}
\begin{proof}
  Let $h\in H^1(\R^3)$ be a stationary solution. Recall that the perturbation $\delta=u-h$ satisfies
  \begin{equation}\label{eq:delta-relative-proof}
    \delta_{tt}+\gamma\delta_t-\Delta\delta+m\delta+V\delta=\mathcal N(\delta),
  \end{equation}
  where
  $$
    \mathcal N(\delta):=-\delta^2\Gamma(\delta)+f_2(h)-f_2(h+\delta) \quad \text{and}\quad \Gamma(\delta)=\int_0^1 \int_0^\tau f^{\prime\prime}_1 ( h+\theta  \delta) \td \theta\td\tau.
  $$

  By the regularity estimate for $\delta$ and the energy bound, we have
  $$
    \delta\in C([0,T];H^1(\R^3)),
    \qquad
    \delta_t\in C([0,T];L^2(\R^3)),
  $$
  and
  $$
    \delta\in
    L^{\frac{8}{p-1}}
    \bigl(0,T;L^{\frac{24}{5-p}}(\R^3)\bigr).
  $$
  Interpolating this space-time bound with $\delta\in L^\infty(0,T;L^6(\R^3))$, which follows from the energy bound, and using $T<\infty$, we obtain
  $$
    \delta\in
    L^2\bigl(0,T;L^{\frac{12}{5-p}}(\R^3)\bigr)
    \cap
    L^p\bigl(0,T;L^{2p}(\R^3)\bigr).
  $$

  The growth assumptions on $f_1$ and the global Lipschitz property of $f_2$
  imply
  $$
    |\mathcal N(z)|
    \le C\bigl(
    |z|+|z|^2|h|^{p-2}+|z|^p
    \bigr)
  $$
  and
  \begin{equation}\label{eq:N-difference-proof}
    \begin{aligned}
      |\mathcal N(z)-\mathcal N(w)|
      \le C|z-w|\bigl(
      1+|h|^{p-2}(|z|+|w|)
      +|z|^{p-1}+|w|^{p-1}
      \bigr).
    \end{aligned}
  \end{equation}
  Thus,
  \begin{equation}\label{eq:N-L1L2-proof}
    \mathcal N(\delta)\in L^1(0,T;L^2(\R^3)).
  \end{equation}
  Moreover, if $z$ and $w$ are bounded in
  $$
    L^\infty(0,T;L^2)
    \cap
    L^2\bigl(0,T;L^{\frac{12}{5-p}}\bigr)
    \cap
    L^p(0,T;L^{2p}),
  $$
  then \eqref{eq:N-difference-proof} and H\"older's inequality yield
  \begin{align}
    \|\mathcal N(z)-\mathcal N(w)\|_{L^1(0,T;L^2)}
    \le C\Bigl(\|z-w\|_{L^\infty(0,T;L^2)}+\|z-w\|_{L^2(0,T;L^{\frac{12}{5-p}})}+\|z-w\|_{L^p(0,T;L^{2p})}
    \Bigr),
    \label{eq:N-continuity-proof}
  \end{align}
  where the constant depends on $T$ and the corresponding bounds of $z$ and $w$ and on $h$.

  We now regularise only in time. Let $\eta\in C_c^\infty(-1,1)$ be nonnegative with
  $\int_{\R}\eta=1$, set $\eta_\eps(t):=\eps^{-1}\eta(t/\eps)$, and define
  $$
    \delta^\eps:=\int_0^T\eta_\eps(t-\tau)\delta(\tau)\td \tau.
  $$
  Then $\delta^\eps$ is smooth in time with values in $H^1(\R^3)$.
  Fix a compact interval $[a,b]\subset(0,T)$ and let $0<\varepsilon<\min\{a,T-b\}$. Then, for every $t\in[a,b]$,
  $$
    \eta_\varepsilon(t)=\eta_\varepsilon(t-T)=0.
  $$
  Hence, for every $t\in[a,b]$, integrating by parts gives
  $$\pt\delta^\eps(t) =\int_0^T\pt\eta_\eps(t-\tau)\delta(\tau)\td \tau=\int_0^T\eta_\eps(t-\tau)\delta_t(\tau)\td \tau-\eta_\eps(t-\tau)\delta(\tau)\big|^T_0=\delta_t^\eps(t).$$
  Then $\delta_t^\eps$ belongs to $H^1(\R^3)$ and
  $$
    \delta_t^\eps\to\delta_t
    \quad\text{in }C([a,b];L^2(\R^3)).
  $$
  Furthermore,
  $$
    \delta^\eps\to\delta
    \quad\text{in }C([a,b];H^1(\R^3))\cap L^2\bigl(a,b;L^{\frac{12}{5-p}}(\R^3)\bigr)\cap L^p(a,b;L^{2p}(\R^3)).
  $$
  Therefore, \eqref{eq:N-L1L2-proof} and \eqref{eq:N-continuity-proof} imply
  \begin{equation}\label{eq:N-moll-conv-proof}
    (\mathcal N(\delta))^\eps\to\mathcal N(\delta) \quad\text{and}\quad
    \mathcal N(\delta^\eps)\to\mathcal N(\delta)
    \quad\text{in }L^1(a,b;L^2(\R^3)).
  \end{equation}

  Convolving \eqref{eq:delta-relative-proof} in time gives
  $$
    \delta_{tt}^\eps+\gamma\delta_t^\eps
    -\Delta\delta^\eps+m\delta^\eps
    +V\delta^\eps
    =
    (\mathcal N(\delta))^\eps
  $$
  in $H^{-1}(\R^3)$ for $t\in[a,b]$.

  Since $\rho,\nabla\rho\in L^\infty(\R^3)$, multiplication by $\rho$
  maps $H^1(\R^3)$ continuously into itself. Hence, $\rho\delta_t^\eps\in H^1(\R^3)$ is an admissible test function. Testing the regularised equation by $\rho\delta_t^\eps$ and integrating
  over $[a,b]$ gives
  \begin{equation}    \label{eq:regularized-relative-energy-proof}
    \begin{aligned}
       & \left[\frac12\int_{\R^3}\rho\left( |\delta_t^\eps|^2 +|\nabla\delta^\eps|^2 +m|\delta^\eps|^2 +V|\delta^\eps|^2
      \right)\,\tdx \right]_a^b                                                                                          \\
       & \quad
      +\gamma\int_a^b\int_{\R^3}
      \rho|\delta_\tau^\eps|^2\,\tdx\,\td\tau
      +\int_a^b\int_{\R^3}
      (\nabla\delta^\eps\cdot\nabla\rho)
      \delta_\tau^\eps\,\tdx\,\td\tau
      =
      \int_a^b\int_{\R^3}
      \rho(\mathcal N(\delta))^\eps
      \delta_\tau^\eps\,\tdx\,\td\tau .
    \end{aligned}
  \end{equation}

  Define
  $$
    \Phi(x,z):=
    F(h(x)+z)-F(h(x))-f(h(x))z-\frac12V(x)|z|^2.
  $$
  Then
  $$
    \partial_z\Phi(x,z)
    =f(h+z)-f(h)-V z
    =-\mathcal N(z).
  $$
  Since $\delta^\eps$ is $C^1$ in time with values in $H^1$, the usual difference-quotient argument gives
  $$
    \frac{\td}{\td\tau}\int_{\R^3}\rho\Phi(x,\delta^\eps(\tau))\,\tdx=-\int_{\R^3}\rho\mathcal N(\delta^\eps)\delta_\tau^\eps\,\tdx
  $$
  for almost every $\tau\in(a,b)$.

  Adding this chain rule to \eqref{eq:regularized-relative-energy-proof}, we obtain
  \begin{equation}\label{eq:combined-relative-energy-proof}
    \bigl[\widetilde E_\rho(\delta^\eps)\bigr]_a^b +\int_a^b\int_{\R^3}\l(\gamma\rho|\delta_\tau^\eps|^2+
    \nabla\delta^\eps\cdot\nabla\rho\delta_\tau^\eps\r)\tdx \td\tau =
    \int_a^b\int_{\R^3} \rho\Bigl((\mathcal N(\delta))^\eps-\mathcal N(\delta^\eps) \Bigr) \delta_\tau^\eps\,\tdx\,\td\tau ,
  \end{equation}
  where
  $$
    \widetilde E_\rho(\delta)
    :=
    \frac12\int_{\R^3}\rho
    \bigl(
    |\delta_t|^2+|\nabla\delta|^2+m|\delta|^2
    \bigr)\,\tdx
    +\int_{\R^3}\rho
    \bigl(
    F(h+\delta)-F(h)-f(h)\delta
    \bigr)\,\tdx .
  $$
  By \eqref{eq:N-moll-conv-proof} and the uniform bound on $\delta_t^\eps$ in $L^\infty(a,b;L^2(\R^3))$, the right-hand side of \eqref{eq:combined-relative-energy-proof} tends to zero.

  We now pass to the limit in the remaining terms. By the strong convergences established above and the growth
  assumptions on $f$, we have
  $$
    \widetilde E_\rho(\delta^\eps(\cdot))
    \to
    \widetilde E_\rho(\delta(\cdot))\quad \text{in}\quad C([a,b]).
  $$
  Moreover,
  $$
    \int_a^b\int_{\R^3}\rho|\delta_\tau^\eps|^2\,\tdx\,\td\tau\to\int_a^b\int_{\R^3}\rho|\delta_\tau|^2\,\tdx\,\td\tau
  $$
  and
  $$
    \int_a^b\int_{\R^3}(\nabla\delta^\eps\cdot\nabla\rho)\delta_\tau^\eps\,\tdx\,\td\tau\to\int_a^b\int_{\R^3} (\nabla\delta\cdot\nabla\rho)\delta_\tau\,\tdx\,\td\tau .
  $$
  Therefore,
  \begin{equation}\label{eq:relative-energy-final-simplified}
    \begin{aligned}
      \widetilde E_\rho(\delta(b))-\widetilde E_\rho(\delta(a))+\gamma\int_a^b\int_{\R^3} \rho|\delta_\tau|^2\,\tdx\,\td\tau
      +\int_a^b\int_{\R^3} (\nabla\delta\cdot\nabla\rho)\delta_\tau\,\tdx\,\td\tau =0 .
    \end{aligned}
  \end{equation}

  It remains to rewrite the relative energy in terms of $E_\rho(u)$. Since
  $\rho u(t)\in H^1(\R^3)$, testing the stationary equation for $h$
  by $\rho u(t)$ yields
  $$
    \begin{aligned}
      \langle g,\rho u(t)\rangle
      =
      \int_{\R^3}
      \bigl(
      \rho\nabla h\cdot\nabla u(t)
      +u(t)\nabla h\cdot\nabla\rho
      +m\rho h u(t)
      +\rho f(h)u(t)
      \bigr)\,\tdx .
    \end{aligned}
  $$
  Using $u=h+\delta$, expanding the squares, and applying this identity, we obtain
  $$
    \begin{aligned}
      \widetilde E_\rho(\delta(t))
      =
       & E_\rho(u(t))
      +\int_{\R^3}
      u(t)\nabla h\cdot\nabla\rho\,\tdx
      +\int_{\R^3}
      \rho\bigl(\frac12|\nabla h|^2+\frac12mh^2+f(h)h-F(h)\bigr)\,\tdx.
    \end{aligned}
  $$
  Hence,
  \begin{equation}\label{eq:weighted-energy-interior}
    E_\rho(u(b))-E_\rho(u(a))+\gamma\int_a^b\int_{\R^3}\rho|u_\tau|^2\,dx\,d\tau+\int_a^b\int_{\R^3}(\nabla u\cdot\nabla\rho)u_\tau\,dx\,d\tau=0.
  \end{equation}

  Finally, the growth assumptions imply
  $$
    |F(s)|\le C\bigl(|s|^2+|s|^{p+1}\bigr),
    \qquad p+1<6.
  $$
  Hence, $E_\rho(u(t))$ is finite and continuous in $t$, because
  $u\in C([0,T];H^1)$ and $\rho\in L^\infty$. Also,
  $$
    \rho|u_t|^2\in L^1((0,T)\times\R^3)
  $$
  and
  $$
    |(\nabla u\cdot\nabla\rho)u_t|
    \le
    \|\nabla\rho\|_{L^\infty}|\nabla u||u_t|
    \in L^1((0,T)\times\R^3).
  $$
  Thus, \eqref{eq:weighted-energy-interior} extends by continuity to every
  $0\le s\le t\le T$:
  $$
    E_\rho(u(t))-E_\rho(u(s))+\gamma\int_s^t\int_{\R^3}\rho|u_\tau|^2\,dx\,d\tau+\int_s^t\int_{\R^3}(\nabla u\cdot\nabla\rho)u_\tau\,dx\,d\tau=0.
  $$

  Since the integrands in the two integral terms belong to $L^1(0,T)$, the function $E_\rho(u)$ is absolutely continuous on $[0,T]$. Differentiating the
  integrated identity yields
  $$
    \frac{\td}{\tdt}E_\rho(u)
    +\gamma\int_{\R^3}\rho|u_t|^2\,\tdx
    +\int_{\R^3}
    (\nabla u\cdot\nabla\rho)u_t\,\tdx
    =0
  $$
  for almost every $t\in[0,T]$.
\end{proof}
\begin{remark}\label{rem:wgh_engy_eq}
  The weighted energy identity provides a rigorous basis for the energy estimates used in the dissipativity lemma and the tail estimate. In the
  dissipativity argument, we take $\rho\equiv 1$; in the tail estimate, we take $\rho=\chi_n$ and combine the weighted identity with the weak
  formulation tested against $\chi_n u$. These arguments yield the corresponding energy inequalities for weak solutions. Thus, the formal
  testing of the equation by $u_t$ in the subsequent estimates is understood as shorthand for these rigorous arguments.
\end{remark}
\subsection{Dissipativity}

\begin{lemma}[Dissipativity]\label{le_disp}
  Let $g\in H^{-1}(\R^3)$ and let the nonlinearity $f$ satisfy \eqref{F1} and \eqref{F2}. Then the semigroup $\{S(t)\}_{t\geq 0}$ possesses a bounded, positively invariant absorbing set $\mathscr{B}_0$ in $\H$: for every bounded $B\subset \H$ there is a $T_0$, depending only on the $\H$-bound of $B$, with
  \begin{eqnarray*}
    S(t)B\subset \mathscr{B}_0 \quad\text{for all } t\ge T_0.
  \end{eqnarray*}
\end{lemma}
\begin{proof}
  By Remark~\ref{rem:wgh_engy_eq}, the formal testing by $u_t$ can be justified by the weighted energy identity and yields the same energy
  identity. We therefore formally test \eqref{eq:problem} by $u_t+\alpha u$ and obtain
  \begin{equation}\label{eq:ut+alu}
    \begin{aligned}
       & \frac{\td}{\td t}E_\al(u)+(\ga-\al)\|u_t\|^2+\al\|\nabla u\|^2+\al m\|u\|^2 +\al\int_{\R^3} f(u)u\tdx-\al\la g,u\ra=0,
    \end{aligned}
  \end{equation}
  where
  $$E_\al(u)=E(u)+\frac{\al\ga}{2}\|u\|^2+\al(u,u_t).$$

  Taking $\al$ small enough, say $0<\al\ll\ga/2$, and using \eqref{F2}, we obtain
  $$
    E_{\al}(u(t))\geq c_1(\|u_t(t)\|^2+\|\nabla u(t)\|^2+\|u\|^2) -c_2(\|g\|_{H^{-1}(\R^3)}^2+1),\qquad t\geq 0,
  $$
  and
  $$
    \begin{aligned}
      (\ga-\al)\|u_t\|^2 & +\al\|\nabla u\|^2+\al m\|u\|^2+\al\int_{\R^3} f(u)u\tdx-\al\la g,u\ra                         \\
                         & \geq c_3(\|u_t(t)\|^2+\|\nabla u(t)\|^2+\|u\|^2)-c_4(\|g\|_{H^{-1}(\R^3)}^2+1),\qquad t\geq 0,
    \end{aligned}
  $$
  where the positive constants $c_i$, $i=1,2,3,4$, depend only on $\gamma$, $\alpha$ and $m$. Here $\alpha>0$ is fixed sufficiently small in terms of $\gamma$ and $m$. The nonlinear terms are used only through $F(u)\geq 0$ and $f(u)u\geq 0$. Integrating \eqref{eq:ut+alu} over $[0,t]$ therefore gives
  \begin{equation}\label{eq:Eal_int}
    \begin{aligned}
      c_1(\|u_t(t)\|^2 & +\|\nabla u(t)\|^2+\|u\|^2) -c_2(\|g\|_{H^{-1}(\R^3)}^2+1)-E_{\al}(u(0))                        \\
                       & \leq-\int_0^t [c_3(\|u_t(s)\|^2+\|\nabla u(s)\|^2+\|u\|^2)-c_4(\|g\|_{H^{-1}(\R^3)}^2+1)]\td s.
    \end{aligned}
  \end{equation}

  Let $M=c_4(\|g\|_{H^{-1}(\R^3)}^2+1)/c_3+1$. For a bounded set $B\subset\H$, the quantity
  $E_\alpha(u(0))$ is uniformly bounded for $(u_0,u_1)\in B$.
  If for every $T>0$ there were an orbit with
  $$
    \|u_t(t)\|^2+\|\nabla u(t)\|^2+\|u(t)\|^2>M
    \qquad\text{for all }t\in[0,T],
  $$
  then evaluating \eqref{eq:Eal_int} at $t=T$ would give a left-hand side bounded below uniformly in $T$, while the right-hand side is
  $\le -c_3T$. This is impossible for large $T$. Hence, there exists $T_B=T(\|B\|_{\H},M)$, independent of the individual orbit, and
  for every $(u_0,u_1)\in B$ some $t_0\in[0,T_B]$ with
  $$
    \|u_t(t_0)\|^2+\|\nabla u(t_0)\|^2+\|u(t_0)\|^2\leq M.
  $$
  The energy estimate \eqref{eq:eng_est} and the autonomy of the equation then give
  $$
    \|(u(t),u_t(t))\|_{\H} \le Q(M,\|g\|_{H^{-1}(\R^3)}) \qquad\text{for all  }t\ge T_B.
  $$
  Thus,
  $$
    \mathscr B=\{(u,v)\in\H:\|(u,v)\|_{\H} \le Q(M,\|g\|_{H^{-1}(\R^3)})\}
  $$
  is a bounded absorbing set.

  Choose $T_*>0$ such that
  $$
    S(t)\mathscr B\subset\mathscr B\qquad\text{for all }t\geq T_*.
  $$
  Let
  $$
    \mathscr B_0=\overline{\bigcup_{\tau\geq T_*}S(\tau)\mathscr B}^{\,\H}.
  $$
  Then $\mathscr B_0$ is bounded, absorbing, and positively invariant:
  \begin{equation}\label{eq:B0_positive_invariant}
    S(t)\mathscr B_0\subset\mathscr B_0 \qquad\text{for all }t\ge 0.
  \end{equation}
\end{proof}

\subsection{Asymptotic compactness}

On an unbounded domain the Sobolev embeddings are no longer compact, so asymptotic compactness cannot be obtained from higher regularity and the noncompact global Sobolev embedding alone. We circumvent this by the uniform tail estimates of \cite{Wang1999}, together with a refined decomposition $u=h+v+w$ in which $v$
decays exponentially while the higher-regularity component $w$ remains bounded in $\H^{1+\vartheta}$.

More precisely, we first show that solutions of \eqref{eq:problem} are uniformly small outside a large ball, and then invoke the Rellich theorem on that ball. It is the energy identity of Lemma \ref{le_energy_eqv} that makes both the dissipativity and the tail estimate available for weak solutions.

For the reader's convenience we recall the definition of the global attractor; see \cite{avb1,rt} for details.
\begin{Def}
  Let $\{S(t)\}_{t\geq 0}$ be a semigroup on a metric space $(E,d)$.
  A subset $\mathscr{A}\subset E$ is a global attractor for the
  semigroup if it is compact and
  \begin{itemize}
    \item[(1)] invariant, that is, $S(t)\mathscr{A}=\mathscr{A}$ for all $t\geq 0$;
    \item[(2)] attracting for all bounded sets of $E$: for every bounded
      $B\subset E$,
      $$ d(S(t)B,\mathscr{A})\rightarrow 0 \quad\text{as } t\rightarrow +\infty,$$
  \end{itemize}
  where $d(B,A)=\sup _{x\in B}\inf_{y\in A}d(x,y)$ denotes the semidistance between two sets $A$ and $B$.
\end{Def}
The following uniform estimate on the tails of the solutions is the key to asymptotic compactness.
\begin{thm}\label{thm_tail}
  Let $g\in H^{-1}(\R^3)$ and let $f\in C^3(\R)$ satisfy \eqref{F1}--\eqref{F3}. Then for every bounded set $B\subset \H$ and every $\eps>0$ there exist $T_0=T_0(\eps, B)$ and $N=N(\eps, B)$ such that every solution $u$ of \eqref{eq:problem} with $(u_0,u_1)\in B$ satisfies
  \begin{equation}
    \int_{|x|\ge n}\l(|u_t(x,t)|^2+|u(x,t)|^2+|\nabla u(x,t)|^2\r)\td x<\eps
  \end{equation}
  for all $t\ge T_0$ and all $n\ge N$.
\end{thm}
\begin{proof}
  Let $\chi: \R^3\to \R$ be smooth with $0\le \chi\le 1$ and
  $$\chi(x)=0 \quad \text{for }|x|\le 1/2\quad\text{and}\quad \chi(x)=1\quad \text{for } |x|\ge 1,$$
  and put $\chi_n(x)=\chi(x/n)$ for $n\in \mathbb{N}$, so that
  $$|\nabla \chi_n(x)|\le C/n$$
  with $C$ independent of $n$.

  Similarly, using the weighted energy identity with $\rho=\chi_n$, we may legitimately employ the formal test function $\chi_n u_t+\alpha\chi_n u$ in \eqref{eq:problem}, which gives
  \begin{equation}
    \begin{aligned}
       & \frac{\td}{\td t}E_{\al,n}+\al \beta E_{\al,n}+(\ga- \al-\frac{\al\beta}{2})\int_{\R^3} \chi_n(x) |u_t(x)|^2 \td x                        \\
       & +(\al-\frac{\al\beta}{2})\int_{\R^3} \chi_n(x)( |\nabla u(x)|^2+m|u(x)|^2)\td x-\frac{\al^2\beta\ga}2\int_{\R^3} \chi_n(x) |u(x)|^2 \td x \\
       & -\al^2\beta\int_{\R^3} \chi_n(x) u_t(x)u(x)\td x +\int_{\R^3}\nabla u(x)\cdot\nabla \chi_n(x) (u_t(x)+\al u(x))\td x                      \\
       & +\al\int_{\R^3} \chi_n(x) f(u)u\tdx-\al\beta\int_{\R^3} \chi_n(x)F(u(x))\td x-\al(1-\beta)\la g,\chi_n u\ra=0,
    \end{aligned}
  \end{equation}
  where the modified energy functional is given by
  \begin{equation*}
    E_{\al,n}=\frac12\int_{\R^3} \chi_n(x)(|u_t(x)+\al u(x)|^2+|\nabla u(x)|^2+(m+\al\ga-\al^2)|u(x)|^2+2F(u(x)))\td x-\la g,\chi_n u\ra.
  \end{equation*}
  Condition \eqref{F3} bounds the nonlinear terms from below,
  \begin{equation}
    \al\int_{\R^3} \chi_n(x) f(u)u\tdx-\al\beta\int_{\R^3} \chi_n(x)F(u(x))\td x\ge 0.
  \end{equation}
  Fix $\al>0$ sufficiently small, depending only on $m$, $\gamma$ and $\beta$, and independent of $n$ and of the initial data. With this fixed $\al$, since $\beta\le 1$ and $\chi_n\ge 0$,
  \begin{equation}
    (\al-\frac{\al\beta}{2})\int_{\R^3} \chi_n(x)( |\nabla u(x)|^2+m|u(x)|^2)\td x-\frac{\al^2\beta\ga}2\int_{\R^3} \chi_n(x) |u(x)|^2 \td x\ge 0.
  \end{equation}
  Moreover, by the Cauchy--Schwarz inequality,
  $$
    \al^2\beta\left|\int_{\R^3} \chi_n u_t u\td x\right|
    \le\frac{\al^2\beta}{2}\int_{\R^3}\chi_n\l(|u_t|^2+|u|^2\r)\td x,
  $$
  which for this fixed $\al$ is absorbed by the two preceding terms. We are left with the differential inequality
  \begin{equation}
    \frac{\td}{\td t}E_{\al,n}+\al \beta E_{\al,n}\le -\int_{\R^3}\nabla u(x)\cdot\nabla \chi_n(x) (u_t(x)+\al u(x))\td x+\al(1-\beta)\la g,\chi_n u\ra,
  \end{equation}
  whose right-hand side we now estimate. For the first term,
  \begin{equation}
    \left|\int_{\R^3}\nabla u(x)\cdot\nabla \chi_n(x) (u_t(x)+\al u(x))\td x\right|\lesssim \frac{1}{n} \|\nabla u\|\|u_t+\al u\|.
  \end{equation}
  For the second term, set $g_0=(1-\Delta)^{-1}g$. Since $g\in H^{-1}(\R^3)$, we have $g_0\in H^1(\R^3)$. Integration by parts gives
  \begin{equation}\label{eq:g_chi}
    \begin{aligned}
      \left|\la g,\chi_n u\ra\right| & =\left|\int_{\R^3}\bigl(g_0\chi_nu+\nabla g_0\cdot\nabla(\chi_nu)\bigr)\,\tdx\right|                                               \\
                                     & =\left|\int_{\R^3}(g_0(x)\chi_n u(x)+u(x)\nabla g_0(x)\cdot\nabla \chi_n(x) +\chi_n(x)\nabla g_0(x)\cdot \nabla u(x) )\td x\right| \\
                                     & \le \l(\|g_0\chi_n\| +\|\nabla g_0\cdot\nabla \chi_n\|\r)\|u\|+\|\chi_n\nabla g_0\|\|\nabla u\|                                    \\
                                     & \le \l(\|g_0\chi_n\| +\|\nabla g_0\cdot\nabla \chi_n\|+\|\chi_n\nabla g_0\|\r)\|u\|_{H^1(\R^3)}.
    \end{aligned}
  \end{equation}
  Since $|\nabla\chi_n|\le C/n$ and $g_0 \in H^1(\R^3)$, both $\|g_0\chi_n\|$ and $\|\chi_n\nabla g_0\|$ tend to $0$ as $n\to\infty$. Hence, given $\eta>0$, there is an $N=N(\eta,B)$ such that for all $n\ge N$
  \begin{equation}
    \frac{1}{n} \|\nabla u\|\|u_t+\al u\|+\l(\|g_0\chi_n\| +\|\nabla g_0\cdot\nabla \chi_n\|+\|\chi_n\nabla g_0\|\r)\|u\|_{H^1(\R^3)}\le \eta ,
  \end{equation}
  and therefore
  \begin{equation}
    \frac{\td}{\td t}E_{\al,n}+\al \beta E_{\al,n}\lesssim \eta.
  \end{equation}
  Since the initial data set $B$ is bounded in $\H$, from the definition of $E_{\al,n}$, the Cauchy--Schwarz inequality, and the boundedness of $B$, we have
  $$
    \sup_{n\in\mathbb{N}} \sup_{(u_0,u_1)\in B} E_{\al,n}(0) \le C(B, \gamma, m, \alpha, g) < \infty,
  $$
  where the nonlinear potential term is controlled by the $H^1$-bound of $u_0$ via the Sobolev embedding.

  Gronwall's inequality then gives
  \begin{equation}
    E_{\al,n}(t)\lesssim C(B, \gamma, m, \alpha, g) e^{-\al\beta t}+ \eta ,
  \end{equation}
  while \eqref{eq:g_chi} yields
  \begin{equation}
    \int_{\R^3} \chi_n(x)(|u_t(x)|^2+|\nabla u(x)|^2+m|u(x)|^2)\td x\lesssim E_{\al,n}(t)+\la g,\chi_n u\ra\le E_{\al,n}(t)+\eta.
  \end{equation}
  With $\al$ already fixed as above, given $\eps>0$, first choose $\eta$ sufficiently small, then choose $N$ so that the cut-off and forcing errors are below $\eta$, and finally choose $T_0$ so that the exponentially decaying term is below $\eta$ for all $t\ge T_0$. Hence,
  $$\int_{\R^3} \chi_n(x)(|u_t(x)|^2+|\nabla u(x)|^2+m|u(x)|^2)\td x\le \eps$$
  for all $t\ge T_0$ and $n\ge N$, with $T_0$ and $N$ large enough. As $\chi_n\equiv1$ on $|x|\ge n$, this is the claimed bound on the exterior energy. Moreover, since $|\nabla\chi_n|\le C/n$,
  $$
    \|\chi_nu\|^2_{H^1(\R^3)}+\|\chi_nu_t\|^2
    \lesssim\int_{\R^3}\chi_n\l(|u_t|^2+|\nabla u|^2+|u|^2\r)\td x+\frac{C}{n^2}\|u\|^2 ,
  $$
  so the same smallness holds for $\|\chi_nu(t)\|_{H^1}+\|\chi_nu_t(t)\|$, which is the form used in the proof of Proposition \ref{prop_omega_comp}.
\end{proof}

Theorem \ref{thm_attr} rests on the following asymptotic compactness property.

\begin{prop}\label{prop_omega_comp}
  Under assumptions of Theorem \ref{thm_attr}, the semigroup $S(t)$ is asymptotically compact on $\H$, i.e., for every bounded sequence $z_j$ in $\H$ and every sequence $t_j\to\infty$, $\{S(t_j)z_j\}$ has a strongly convergent subsequence in $\H$.
\end{prop}

The rest of the section is devoted to its proof, which begins with a refined decomposition of the solution. Let $h\in H^1(\R^3)$ be a stationary solution and put $\delta:=u-h$. Recall from Section~\ref{sec:wp} the splitting $f=f_1+f_2$, with $|f_1^{\prime\prime\prime}(s)|\le C|s|^{p-3}$ and $|f_2^{\prime}(s)|\le C_{\Lambda}$, the potential $V=f_1^\prime(h)$, and the coefficient
\begin{equation}\label{F1_def}
  \Gamma=\int_0^1 \int_0^\tau f^{\prime\prime}_1 ( h+\theta \delta) \td \theta\td\tau .
\end{equation}
Then $\delta$ satisfies
\begin{equation}\label{eq:delta00}
  \delta_{tt} +\ga \delta_t-\Delta \delta+V\delta+m\delta= -\delta^2 \Gamma+f_2(h)-f_2(u).
\end{equation}
For any $(u_0,u_1)\in \mathscr{B}_0$, we decompose $\delta$ as
\begin{equation}\label{eq:u_decomp}
  \delta(t)=v(t)+w(t),
\end{equation}
where the function $v$ solves the linearised equation
\begin{equation}\label{eq:v}
  \begin{cases}
    v_{tt}+ \ga v_{t}-\Delta v+Vv+m v=0, \\
    v(0)=u_0-h,~v_t(0)=u_1,
  \end{cases}
\end{equation}
and the remainder $w$ satisfies
\begin{equation}\label{eq:w}
  \begin{cases}
    w_{tt}+\ga w_{t}-\Delta w+Vw+m w=-\delta^2 \Gamma+f_2(h)-f_2(u), \\
    w(0)=0, ~w_{t}(0)=0.
  \end{cases}
\end{equation}

The component $v$ decays exponentially.
\begin{lemma}\label{le_v_decay}
  There exist $C,\mu>0$, independent of $(u_0,u_1)$ and depending only on $\ga$, $m$ and $\|h\|_{H^1(\R^3)}$, such that the solution $v$ of \eqref{eq:v}
  satisfies
  \begin{equation}\label{eq:v_decay}
    \| v_t(t)\|+\| v(t)\|_{H^1(\R^3)}\le C e^{-\mu t}\l(\| u_0\|_{H^1(\R^3)}+\| u_1\|+\| g\|_{H^{-1}(\R^3)}\r),
    \qquad t\ge0 .
  \end{equation}
\end{lemma}
\begin{proof}
  Pairing \eqref{eq:v} with $v_t+\al v$, $\al>0$, and integrating over $\R^3$ gives
  $$
    \frac{\td}{\td t}\mathcal E_\al(v)+(\ga-\al)\|v_t\|^2
    +\al\l(\|\nabla v\|^2+\int_{\R^3}Vv^2\td x+m\|v\|^2\r)=0,
  $$
  where
  $$\mathcal E_\al(v)=\frac12\|v_t+\al v\|^2+\frac12\l(\|\nabla v\|^2+\int_{\R^3}Vv^2\td x\r)
    +\frac{m+\al(\ga-\al)}{2}\|v\|^2 .$$
  By \eqref{eq:V_coercive} of Lemma \ref{lem_V},
  $$
    \|\nabla v\|^2+\int_{\R^3}Vv^2\td x+m\|v\|^2\ge\frac12\|\nabla v\|^2+m\|v\|^2
    \ge c\,\|v\|^2_{H^1(\R^3)},
  $$
  so that for $0<\al<\ga/2$ small enough $\mathcal E_\al$ is equivalent to
  $\|v_t\|^2+\|v\|^2_{H^1}$ and satisfies
  $\frac{\td}{\td t}\mathcal E_\al(v)+2\mu\,\mathcal E_\al(v)\le0$ for some $\mu>0$. Gronwall's
  inequality, the bound $\|v(0)\|_{H^1}+\|v_t(0)\|\le\|u_0\|_{H^1}+\|u_1\|+\|h\|_{H^1}$ and
  $\|h\|_{H^1}\le Q(\|g\|_{H^{-1}(\R^3)})$ now give \eqref{eq:v_decay}.
\end{proof}

We turn to the $w$-component, whose treatment depends strongly on the growth rate of the nonlinearity.
\begin{lemma}\label{le_w_reglr}
  Let $\vartheta=\frac{10-2p}{p+7}>0$. There exists a constant
  $$
    R=R\bigl(\mathscr B_0,\|g\|_{H^{-1}(\R^3)},m,\gamma,f,\Lambda\bigr)
  $$
  such that
  $$
    \sup_{z\in\mathscr B_0}\sup_{t\ge0}\left(\|w(t)\|_{H^{1+\vartheta}(\R^3)} +\|w_{t}(t)\|_{H^\vartheta(\R^3)}\right)\leq R.
  $$
\end{lemma}
\begin{proof}
  We prove the uniform boundedness of $(w(t),w_t(t))$ in $\H^{1+\vartheta}:=H^{1+\vartheta}(\R^3)\times H^{\vartheta}(\R^3)$ by energy estimates in Besov spaces. Since $\mathscr B_0$ is positively invariant, every trajectory
  starting from $\mathscr B_0$ remains in $\mathscr B_0$. Hence,
  $$
    \sup_{z\in\mathscr B_0}\sup_{t\ge0}\|S(t)z\|_{\H}<\infty.
  $$
  Moreover, Theorem~\ref{regular} gives
  $$
    \sup_{z\in\mathscr B_0}\sup_{t\ge0}\|u_z-h\|_{L^{\frac{8}{p-1}}(t,t+1;L^{\frac{24}{5-p}}(\R^3))} <\infty.
  $$
  Therefore, the quantity
  $$
    R:=1+\|h\|_{H^1(\R^3)}+\sup_{z\in\mathscr B_0}\sup_{t\ge0}\left\{\|S(t)z\|_{\H}+\|u_z-h\|_{L^{\frac{8}{p-1}}(t,t+1;L^{\frac{24}{5-p}}(\R^3))}\right\}
  $$
  is finite and depends only on the fixed data and the absorbing set.

  We claim that the source term $-\delta^2 \Gamma + f_2(h) - f_2(u)$ in \eqref{eq:w} belongs to $L^1_{loc}(\R^{+};\dot{H}^{\vartheta}(\R^3))$. To estimate $\delta^2\Gamma$, we use the fractional Gagliardo--Nirenberg inequality \cite{brezis18}:
  \begin{equation*}
    \|\phi\|_{\dot{H}^{\vartheta}(\R^3)} \lesssim \|\phi\|_{\dot{W}^{1,\frac{8}{p+3}}(\R^3)}^{\vartheta} \|\phi\|_{L^{\frac{12}{p+1}}(\R^3)}^{1-\vartheta},
  \end{equation*}
  where the value $\vartheta=\frac{10-2p}{p+7}$ comes from the interpolation relation $\frac12=\frac{\vartheta(p+3)}{8}+\frac{(1-\vartheta)(p+1)}{12}$. Applying this inequality to $\phi=\delta^2\Gamma$ gives
  \begin{equation*}
    \|\delta^2 \Gamma\|_{\dot{H}^{\vartheta}(\R^3)} \lesssim \|\delta^2 \Gamma\|_{\dot{W}^{1,\frac{8}{p+3}}(\R^3)}^{\vartheta} \|\delta^2 \Gamma\|_{L^{\frac{12}{p+1}}(\R^3)}^{1-\vartheta}.
  \end{equation*}
  We estimate the two factors separately:
  \begin{equation}
    \begin{aligned}
      \|\delta^2 \Gamma\|_{\dot{W}^{1,\frac{8}{p+3}}(\R^3)}
       & = \|\nabla(\delta^2 \Gamma)\|_{L^{\frac{8}{p+3}}(\R^3)}
      = \|2 \delta \Gamma\nabla \delta+\delta^2\nabla \Gamma\|_{L^{\frac{8}{p+3}}(\R^3)}                                                                           \\
       & \le \left(2\|\nabla \delta\|_2 \|\Gamma\|_{L^{\frac{6}{p-2}}}+\|\delta\|_{L^6} \|\nabla \Gamma\|_{L^{\frac{6}{p}}}\right) \|\delta\|_{L^{\frac{24}{5-p}}} \\
       & \lesssim \left(2R \cdot R^{p-2} + R \cdot R^{p-2}\right) \|\delta\|_{L^{\frac{24}{5-p}}}                                                                  \\
       & \lesssim R^{p-1} \|\delta\|_{L^{\frac{24}{5-p}}},
    \end{aligned}
  \end{equation}
  and
  \begin{equation}
    \|\delta^2 \Gamma\|_{L^{\frac{12}{p+1}}(\R^3)} \lesssim \|\delta\|_{L^{\frac{24}{5-p}}(\R^3)}^2 \|\Gamma\|_{L^{\frac{6}{p-2}}(\R^3)}.
  \end{equation}

  Since $\Gamma = \int_0^1\int_0^\tau f_1^{\prime\prime}(h + \theta\delta)\td\theta \td\tau$,
  \begin{equation*}
    \begin{aligned}
      \|\nabla \Gamma\|_{L^{\frac{6}{p}}(\R^3)}
       & \lesssim \sup_{0\le\theta\le 1} \|f_1^{\prime\prime\prime}(h + \theta\delta)(\nabla h + \theta\nabla\delta)\|_{L^{\frac{6}{p}}(\R^3)} \\
       & \lesssim \sup_{0\le\theta\le 1}\|h + \theta\delta\|_{H^1(\R^3)} \||h + \theta\delta|^{p-3}\|_{L^{\frac{6}{p-3}}(\R^3)} \lesssim R^{p-2}.
    \end{aligned}
  \end{equation*}
  Altogether,
  \begin{equation}\label{F1_bound}
    \|\delta^2 \Gamma\|_{\dot{H}^{\vartheta}(\R^3)} \lesssim R^{p-2+\vartheta} \|\delta\|_{L^{\frac{24}{5-p}}(\R^3)}^{2-\vartheta}.
  \end{equation}
  The remaining terms are handled by the embedding $H^1(\R^3) \hookrightarrow \dot{H}^{\vartheta}(\R^3)$:
  \begin{equation*}
    \begin{aligned}
      \|f_2(u)- f_2(h)\|_{H^1(\R^3)}
       & \le  \|\nabla(f_2(u)- f_2(h))\|+\|f_2(u)- f_2(h)\|         \\
       & \le C_\Lambda\bigl(\|\nabla u\|+\|\nabla h\|+\|u-h\|\bigr) \\
       & \le C_\Lambda R\le C_\Lambda R^p.
    \end{aligned}
  \end{equation*}
  Using the a priori bound \eqref{eq:regu_est} on $\delta$ in $L^{\frac{8}{p-1}}(t,t+1;L^{\frac{24}{5-p}}(\R^3))$ gives
  \begin{equation}\label{eq:bd_source}
    \begin{aligned}
      \|\delta^2 \Gamma + f_2(u)- f_2(h)\|_{L^1(t,t+1;\dot{H}^{\vartheta}(\R^3))} & \lesssim C_\Lambda\sup_{t \ge 0} \left(R^{p-2+\vartheta}\|\delta\|_{L^{\frac{8}{p-1}}(t,t+1;L^{\frac{24}{5-p}}(\R^3))}^{2-\vartheta} + R^p\right) \\
                                                                                  & \lesssim C_\Lambda R^{p}.
    \end{aligned}
  \end{equation}
  By the equivalence in Theorem~\ref{eqv_besov}, this implies
  $$
    \sop^{\vartheta/2}(\delta^2 \Gamma + f_2(u)- f_2(h))\in L^1_{\mathrm{loc}}(\R^+;L^2(\R^3)),
  $$
  with
  $$
    \|\sop^{\vartheta/2}(\delta^2 \Gamma + f_2(u)- f_2(h))\|_{L^1(t,t+1;L^2(\R^3))} \lesssim \|\delta^2 \Gamma + f_2(u)- f_2(h)\|_{L^1(t,t+1;\dot H^\vartheta(\R^3))}.
  $$

We next establish the $H^{\vartheta}$-regularity of $(w,w_t)$. We first regularise \eqref{eq:w} by applying the bounded spectral truncation $\sop^{\vartheta/2}\mathbf{1}_{[0,N]}(\sop)$. The resulting energy identity is first established for the regularised equation. The corresponding truncated source terms converge in
$L^1(0,T;L^2(\R^3))$ for every $T>0$, and continuous dependence therefore gives convergence of the associated solutions in $C([0,T];\H)$. Moreover,
$$
\mathbf{1}_{[0,N]}(\sop)(w,w_t) \longrightarrow (w,w_t) \quad\text{in }C([0,T];L^2(\R^3)\times L^2(\R^3)).
$$
By the closedness of $\sop^{\vartheta/2}$, the limiting quantities are
$$
\bigl(\sop^{\vartheta/2}w,\sop^{\vartheta/2}w_t\bigr).
$$
Thus, after passing to the limit $N\to\infty$, we may test \eqref{eq:w} with $\sop^\vartheta(w_t+\al w)$. This gives
  \begin{equation*}
    \begin{aligned}
      \frac{d}{dt}E_{w,\al}(t) + 2\al E_{w,\al}(t) + & (\ga - 2\al)\left\|\sop^{\vartheta/2}w_t + \al\sop^{\vartheta/2} w\right\|^2                             \\
      =                                              & \int_{\R^3} \sop^{\vartheta/2}(-\delta^2 \Gamma + f_2(h) - f_2(u)) \sop^{\vartheta/2}(w_t + \al w)\td x,
    \end{aligned}
  \end{equation*}
  with the modified energy
  \begin{equation*}
    E_{w,\al}(t) = \frac12\left\|\sop^{\vartheta/2}(w_t + \al w)\right\|^2 + \frac12\left\|\sop^{(1+\vartheta)/2}w\right\|^2 + \frac{m - \al(\ga - \al)}{2}\left\|\sop^{\vartheta/2}w\right\|^2.
  \end{equation*}

  For $\al>0$ small enough,
  \begin{equation}\label{eq:egy_}
    \left\|\sop^{\vartheta/2}w_t\right\|^2+\left\|\sop^{(1+\vartheta)/2}w\right\|^2+\left\|\sop^{\vartheta/2}w\right\|^2\lesssim E_{w,\al}(t) .
  \end{equation}
  Choose $\alpha>0$ so that
  $$\gamma-2\alpha\ge0, \qquad m-\alpha(\gamma-\alpha)>0.$$
  Applying the Cauchy--Schwarz inequality, and regularising $\sqrt{E_{w,\al}}$ by $\sqrt{E_{w,\al}+\varepsilon}$ before letting $\varepsilon\downarrow0$, we obtain
  \begin{equation*}
    \frac{\td}{\td t}\sqrt{E_{w,\al}(t)} + \al \sqrt{E_{w,\al}(t)} \le \left\|\sop^{\vartheta/2}(\delta^2 \Gamma + f_2(u)- f_2(h))\right\|.
  \end{equation*}
  As $w(0) = w_t(0) = 0$ we have $E_{w,\al}(0) = 0$, so Gronwall's inequality gives
  \begin{equation}\label{eq:Eal}
    \sqrt{E_{w,\al}(t)} \lesssim \int_0^t e^{\al(s-t)} \left\|\sop^{\vartheta/2}(\delta^2 \Gamma + f_2(u)- f_2(h))\right\|  \td s.
  \end{equation}
  Splitting the integral over unit time intervals, we estimate the right-hand side by
  \begin{equation}
    \begin{aligned}
               & \int_0^{t}e^{\al(s-t)}\left\|\sop^{\vartheta/2}(\delta^2 \Gamma + f_2(u)- f_2(h))\right\|\td s                                         \\
      \le      & \sum_{k=0}^{[t]}\int_k^{k+1}e^{\al(s-t)}\left\|\sop^{\vartheta/2}(\delta^2 \Gamma + f_2(u)- f_2(h))\right\|\td s                       \\
      \le      & \sum_{k=0}^{[t]}e^{\al(k+1-t)}\int_k^{k+1}\left\|\sop^{\vartheta/2}(\delta^2 \Gamma + f_2(u)- f_2(h))\right\|\td s                     \\
      \le      & \sum_{k=0}^{[t]}e^{\al(k+1-t)} \sup_{k\ge 0}\left\|\sop^{\vartheta/2}(\delta^2 \Gamma + f_2(u)- f_2(h))\right\|_{L^1(k,k+1;L^2(\R^3))} \\
      \lesssim & \sup_{k\ge 0}\left\|\sop^{\vartheta/2}(\delta^2 \Gamma + f_2(u)- f_2(h))\right\|_{L^1(k,k+1;L^2(\R^3))}.
    \end{aligned}
  \end{equation}
  Substituting this into \eqref{eq:Eal} gives
  \begin{equation}
    \left\|\sop^{\vartheta/2}w_t\right\|+\left\|\sop^{(1+\vartheta)/2}w\right\|\lesssim \sup_{k\ge 0}\left\|\sop^{\vartheta/2}(\delta^2 \Gamma + f_2(u)- f_2(h))\right\|_{L^1(k,k+1;L^2(\R^3))}.
  \end{equation}

  Since $\vartheta>0$ and $1+\vartheta<\frac32$, Theorem~\ref{eqv_besov} implies the norm equivalences
  $$
    \dot{B}^{0}_{2,2}(\sqsop)\cong L^2(\R^3),\qquad \dot{B}^{\vartheta}_{2,2}(\sqsop)\cong \dot{H}^{\vartheta}(\R^3),\qquad \dot{B}^{1+\vartheta}_{2,2}(\sqsop)\cong \dot{H}^{1+\vartheta}(\R^3).
  $$
  Moreover, the spectral calculus gives
  $$
    \sop^{\vartheta/2}: \dot{B}^{\vartheta}_{2,2}(\sqsop)\longrightarrow \dot{B}^{0}_{2,2}(\sqsop),
    \qquad
    \sop^{(1+\vartheta)/2}: \dot{B}^{1+\vartheta}_{2,2}(\sqsop)\longrightarrow \dot{B}^{0}_{2,2}(\sqsop),
  $$
  with equivalent norms. Applying these equivalences to the preceding spectral energy estimate, we obtain
  \begin{equation}\label{eq:key_estimate}
    \|w_t(t)\|_{\dot{H}^{\vartheta}(\R^3)} +\|w(t)\|_{\dot{H}^{1+\vartheta}(\R^3)}
    \lesssim   \sup_{k\ge0} \|\delta^2\Gamma+f_2(u)-f_2(h)\|_{L^1(k,k+1;\dot{H}^{\vartheta}(\R^3))}.
  \end{equation}
  Combining \eqref{eq:key_estimate} with \eqref{eq:bd_source} yields
  $$
    \|w_t(t)\|_{\dot{H}^{\vartheta}(\R^3)}+\|w(t)\|_{\dot{H}^{1+\vartheta}(\R^3)}\lesssim C_\Lambda R^p.
  $$
  Together with the uniform $\H$-bound $\|w(t)\|_{H^1(\R^3)}+\|w_t(t)\|\lesssim R$, a consequence of \eqref{eq:u_decomp}, \eqref{eq:v_decay} and Lemma \ref{le_disp}, this controls the non-homogeneous norms,
  $$\|w_t(t)\|_{H^{\vartheta}(\R^3)} + \|w(t)\|_{H^{1+\vartheta}(\R^3)}\lesssim C_\Lambda R^{p}.$$

  This exhibits a bounded set $\mathcal{C} \subset \H^{1+\vartheta}(\R^3)$ containing $(w(t), w_t(t))$ for all $t \ge 0$.  The set $\mathcal{C}$ depends on $\Lambda$, hence on $\|g\|_{H^{-1}(\R^3)}$.
\end{proof}
Let $\chi^c_n = 1 - \chi_n$, where $\chi_n$ is the cut-off function introduced in the proof of Theorem \ref{thm_tail}. Since multiplication by a fixed $C_c^\infty$ cut-off is bounded on $H^s(\R^3)$ for $s\ge0$, the localised component remains bounded in $\H^{1+\vartheta}(\R^3)$.
\begin{lemma}\label{le_chiw_reglr}
  For every $n\in \N$ there is a bounded set $\mathcal{C}_n \subset \H^{1+\vartheta}(\R^3)$ with $(\chi^c_n w(t), \chi^c_n w_t(t))\in \mathcal{C}_n$ for all $t \ge 0$.
\end{lemma}
\begin{proof}
  Since $\chi_n(x)=\chi(x/n)$ with $\chi = 0$ on $|x|\le 1/2$ and $\chi = 1$ on $|x|\ge 1$, the
  function $\chi^c_n = 1-\chi_n$ belongs to $C_c^\infty(\R^3)$, equals $1$ on $B_{n/2}(0)$
  and vanishes outside $B_n(0)$. Multiplication by a fixed function of
  $C_c^\infty(\R^3)$ is bounded on $H^{s}(\R^3)$ for every $s\ge0$ (see, e.g.,
  \cite{bcd}), so that for some $C_n>0$
  $$
    \|\chi^c_n w\|_{H^{1+\vartheta}(\R^3)} \le C_n \|w\|_{H^{1+\vartheta}(\R^3)}
    \quad\text{and}\quad
    \|\chi^c_n w_t\|_{H^{\vartheta}(\R^3)} \le C_n \|w_t\|_{H^{\vartheta}(\R^3)} .
  $$
  Combined with Lemma \ref{le_w_reglr}, these bounds place $(\chi^c_n w(t), \chi^c_n w_t(t))$
  in a bounded subset $\mathcal{C}_n$ of $\H^{1+\vartheta}$ for all $t \ge 0$.
\end{proof}

We can now complete the proof of Proposition~\ref{prop_omega_comp}. Let $\eps>0$. Choose $n=n(\eps)$ and $T=T(\eps)$
so large that, for all $z=(u_0,u_1)\in\mathscr B_0$ and all $t\ge T$,
$$
  \|(\chi_n u_z(t),\chi_n u_{z,t}(t))\|_{\H}\le C\eps,
  \qquad
  \|(v_z(t),v_{z,t}(t))\|_{\H}\le C\eps,
$$
where $u_z$ denotes the solution with initial datum $z$, and $u_z=h+v_z+w_z$. This is possible by Theorem~\ref{thm_tail} and
Lemma~\ref{le_v_decay}.

Set $\chi_n^c=1-\chi_n$ and define
$$
  K_n:=\overline{
    \{(\chi_n^c(h+w_z(t)),\chi_n^c w_{z,t}(t)):
    z\in\mathscr B_0,\ t\ge0\}}
  ^{\,\H}.
$$
By Lemmas~\ref{le_w_reglr} and~\ref{le_chiw_reglr}, $(\chi_n^c w_z,\chi_n^c w_{z,t})$ is bounded in
$\H^{1+\vartheta}$ and supported in $B_n(0)$. By the Rellich theorem and the fact that adding the fixed vector $(\chi_n^c h,0)\in\H$ preserves compactness, $K_n$ is compact in $\H$.

For $t\ge T$ and $z\in\mathscr B_0$,
$$
  S(t)z =(\chi_n^c(h+w_z),\chi_n^c w_{z,t}) +(\chi_n^c v_z,\chi_n^c v_{z,t}) +(\chi_n u_z,\chi_n u_{z,t}),
$$
so
$$
  \sup_{z\in\mathscr B_0,\ t\ge T} \operatorname{dist}_{\H}\l(S(t)z,K_n\r) \le C\eps .
$$
Since $K_n$ is compact, it has a finite $C\eps$-net. Consequently, the whole tail $\bigcup_{t\ge T}S(t)\mathscr B_0$ is covered by finitely many
$C'\eps$-balls in $\H$.

Now let $\{z_j\}\subset\H$ be bounded and $t_j\to\infty$. Set $B=\{z_j:j\in\N\}$. By Lemma~\ref{le_disp}, there is a
fixed $T_B\ge0$ such that $S(T_B)B\subset\mathscr B_0$. After discarding finitely many terms so that $t_j\ge T_B$, put
$$  y_j=S(T_B)z_j\in\mathscr B_0,\qquad \tau_j=t_j-T_B\to\infty.$$
By the semigroup property, $S(t_j)z_j=S(\tau_j)y_j$. Applying the above finite-covering property with $\eps=2^{-k}$, $k=1,2,\dots$, we extract a diagonal subsequence, which we do not relabel, such that $\{S(\tau_j)y_j\}$ is Cauchy in $\H$.
Since $\H$ is complete, this subsequence converges strongly in $\H$. We have thus proved that $S(t)$ is asymptotically compact in $\H$.

\begin{proof}[Proof of Theorem \ref{thm_attr}]
  Lemma~\ref{le_disp} supplies a bounded, positively invariant absorbing set $\mathscr B_0$ for the semigroup $\{S(t)\}_{t\ge0}$.
  Theorem~\ref{thm_uniq} ensures that each map $S(t):\H\to\H$ is continuous, and Proposition~\ref{prop_omega_comp} establishes
  asymptotic compactness in $\H$. The standard abstract attractor theorem for continuous semigroups with a bounded absorbing set
  and asymptotic compactness (see, e.g., \cite{rt}) therefore implies that
  $$\mathscr A:=\omega(\mathscr B_0)$$
  is compact, invariant, and attracts every bounded subset of $\H$. Hence, $\mathscr A$ is the compact global attractor.
\end{proof}

\section*{Funding}
This work was supported by the NSFC (11801228, 12026431, 12442050, 12671286), the QingLan Project of Jiangsu Province and the Zhongwu Young Innovative Talents Support Program.

\end{document}